\documentclass[12pt]{article}
\usepackage{amssymb,latexsym,amsmath}
\def\C{\mathbb C}
\def\R{\mathbb R}
\def\N{\mathbb N}
\def\Z{\mathbb Z}

\newtheorem{thm}{Theorem}[section]
\newtheorem{lem}{Lemma}[section]

\begin{document}
\sffamily

\title{Zeros of derivatives of strictly\\ non-real meromorphic functions }
\author{J.K. Langley}
\maketitle

\begin{abstract}
 A number of results are proved concerning the existence of non-real zeros of derivatives of strictly non-real meromorphic functions in the plane.
\end{abstract}

\section{Introduction}

Let $f$ be a   meromorphic function in the plane and let $ \widetilde f (z) = \overline{ f( \bar z)} $ 
(this notation will be used throughout). 
Here 
$f$ is called real if $  \widetilde f =f$, and strictly non-real if $\widetilde f $ is not a constant multiple of $f$.
There has been substantial research concerning  non-real zeros of derivatives of real entire or real meromorphic functions
\cite{BEpolya,BEL,CCS2,HSW,HelW1,KiKim,Laams09,rossireal,SS}, 
but somewhat less in the strictly non-real case. The following theorem was proved in \cite{HSWsnr}. 

\begin{thm}[\cite{HSWsnr}]
 \label{thmsnr}
Let $f$ be a strictly non-real meromorphic function in the plane with only real poles.
Then 
$f$, $f'$ and $f''$ have only real zeros if and only if  $f$ has one of the following forms:
\begin{eqnarray*}
(I) \quad f(z) &=& Ae^{Bz}  \quad  ; \\
(II) \quad f(z) &=& A \left(e^{i(cz+d)} - 1 \right)  \quad  ; \\
(III) \quad f(z) &=& A \exp( \exp( i(cz+d) ))  \quad  ; \\
(IV) \quad f(z) &=& A \exp \left[ K( i(cz+d) - \exp( i(cz+d) ) ) \right]  \quad  ; \\
(V) \quad f(z) &=&  \frac{A \exp [ - 2i (cz+d) - 2 \exp( 2i(cz+d)) ] }{\sin^2 (cz+d)}    \quad  ; \\
(VI) \quad f(z) &=&  \frac{A}{e^{i(cz+d)}-1} \quad .
\end{eqnarray*}
Here $A, B \in \C$, while $c, d$ and $K$ are real with $cB \neq 0$ and $K \leq -1/4$. 
\end{thm}

The first aim of the present paper is to prove a result in the spirit of Theorem \ref{thmsnr}, but with no 
assumption on the location of poles.
In \cite{Hin1,Hin2,Hin3} Hinkkanen determined all meromorphic functions  $f$ in the plane such that $f$
and all its  derivatives have
only real zeros, using the fact that under these hypotheses $f$ has at most two distinct poles,
by the P\'olya shire theorem \cite[Theorem 3.6]{Hay2}.
For strictly non-real functions, the following two theorems  will be proved.

\begin{thm}
 \label{thm1}
Let $f$ be a strictly non-real
meromorphic function in the plane such that  all but finitely many zeros of
$f^{(m)}$ are real for $m=0, \ldots, 12$. Then either 
$f'/f$ is a rational function or 
\begin{equation}
 \label{thm1con0}
f(z) = B \left( 1 - T(z) e^{i Az} \right), \quad A \in \R, \quad B \in \C, \quad AB \neq 0, 
\end{equation}
where $T$ is a rational function with $|T(x)| = 1$ for all  $x \in \R$. 

If, in addition, $f$, $f'$, $f''$ and $ f'''$ have only real zeros, then 
$f$ is given by 
one of the following, 
in which $a, b, c , d \in \C$ 
and $\mu \in \Z $: 
\begin{equation}
 \label{thm1con1}
(i)  \quad  f(z) =  (az+b)^{\mu}  \, ; \quad (ii)  \quad f(z) = \frac{az+b}{cz+d}  \,  ; \quad (iii)  \quad f(z) = e^{az+b} - c . 
\end{equation}
\end{thm}

\begin{thm}\label{thm2}
Let $f$ be a strictly non-real
meromorphic function in the plane such that  all  zeros of
$f^{(m)}$ are real for $m=0, \ldots, 9$. 
Then
$f$ is given by (\ref{thm1con1}).
\end{thm}

It is very unlikely that Theorems \ref{thm1} and \ref{thm2}
are sharp in terms of the number of derivatives considered, but examples  (III)-(VI)  of Theorem~\ref{thmsnr} show that 
the absence of non-real zeros of $f$, $f'$ and $f''$  is not enough to imply  (\ref{thm1con1}). 

The next result concerns the   determination of  all meromorphic functions $f$ in the plane  such that $f$ 
and $f''$ have only real zeros and poles (thus  discarding the hypothesis in  Theorem~\ref{thmsnr} that $f'$ has only real zeros).
Such a classification is not known in the real meromorphic case, except when $f$ has finitely many poles \cite{BEL,SS},
or finitely many zeros \cite{Lawiman13},
but for strictly non-real functions of finite lower order the problem is solved by the following theorem,
in which the  terminology is from \cite{Hay2}.


\begin{thm}
 \label{thmsnr17}
Let $f$ be a strictly non-real meromorphic function in the plane, 
such that all but finitely many zeros and poles of $f$ and $f''$ are real. Then $f$ satisfies, as $r \to \infty$,
\begin{equation}
 \label{f'festa}
\overline{N}(r, f) + \overline{N}(r, 1/f)
= O(r) \quad \hbox{and} \quad 
T(r, f'/f) = O( r \log r ) .
\end{equation}
If, in addition, $f$ has finite lower order and all  zeros and poles of $f$ and $f''$ are real
then $f$ is given by one of 
\begin{eqnarray}
 \label{frepa}
(a) \quad f(z) &=& e^{A_1 z + B_1} , \nonumber \\
(b) \quad f(z) &=& e^{ia_1 z} ( T_1'(z) \sin (a_1 z + b_1) - T_1(z) \cos (a_1 z + b_1) ) , \\
(c) \quad f(z) &=& \frac{ T_1(z) }{e^{2i (a_1z+b_1)} - 1 }, \nonumber 
\end{eqnarray}
in which $A_1, B_1 \in \C$ and $a_1, b_1 \in \R$, while $T_1$ is a polynomial of degree at most $1$
such that $T_1(z) = 0$ implies $\sin (a_1z+b_1)=0$. 
\end{thm}

If $T_1$ is a non-zero constant in  (b) or (c) of (\ref{frepa}) then $f$ reduces to (II) or (VI) of Theorem~\ref{thmsnr}
and $f$, $f'$ and $f''$ all have only real zeros and poles. 
However, $T_1$ is non-constant in both of the following examples:
$$
f_1(z) = e^{iz} ( \sin z - z \cos z ), \quad f_1'(z) = e^{iz} ( z \sin z + i ( \sin z - z \cos z )), \quad f_1''(z) = 2z e^{2i z}  ;
$$
$$
f_2(z) = \frac{z}{e^{2iz}-1} , \quad f_2'(z) = \frac{(1-2iz)e^{2iz} - 1}{(e^{2iz}-1)^2} ,
\quad f_2''(z) = \frac{(4i-4z)e^{2iz} - (4i+4z)e^{4iz}}{(e^{2iz}-1)^3} .
$$
Here $f_1'$ and $f_2'$ each have infinitely many non-real zeros, but   $f_2''(z) = 0$ forces 
$$
e^{2iz} = \frac{4i-4z}{4i+4z} , \quad z = \tan z ,
$$
all solutions of which are real (see Lemma \ref{lemtan}), as are all zeros of $f_1$. Furthermore, writing
$$
\frac{f_3'(z)}{f_3(z)} = \frac1{i+z} + i e^{iz} ,
\quad \frac{f_3''(z)}{f_3(z)} = \left( \frac{i-z}{i+z} \right) e^{iz} - e^{2iz} ,
$$
defines a strictly non-real entire function $f_3$ of infinite order, with one zero,
for which $f_3'$ has infinitely many non-real zeros, while all but finitely many
zeros of $f_3''$ are real by \cite[Lemma 2.3]{LaJdA2013}. 

The author thanks John Rossi for  helpful discussions, and the referee for a very careful reading of the manuscript and several helpful suggestions.

\section{Preliminaries}

The following theorem is a combination of results from \cite{FHP,FL2,La5} and uses notation from \cite{Hay2}. 

\begin{thm}
[\cite{FHP,FL2,La5}]
Let $h$ be a non-constant meromorphic function in the plane. \\
(i) 
For $n \geq 3$ there exists $c_n > 0$, depending only on $n$,  such that  
$$
T(r, h'/h) \leq c_n \left( \overline{N}(r, 1/h) + \overline{N}\left(r, 1/h^{(n)}\right) \right) + O( \log r )  \quad \hbox{ as $r \to \infty$. }
$$
(ii) If $n \geq 2$ and $h$ and $h^{(n)}$ have finitely many zeros, then $h'/h$ is a rational function: equivalently,  $h=Se^Q$ with 
$S$ a rational function and $Q$ a polynomial. 
\label{FHPthm}
\end{thm}
Here part (i) follows from \cite[Theorem 3]{FL2} (which should be stated for functions which have
transcendental logarithmic derivative, rather than merely being themselves transcendental), and part (ii) was proved in \cite{FHP,La5}.

\begin{thm}[\cite{BLa}]\label{norfam}
Let $k \geq 2$ and let ${\mathcal H}$ be a family of functions meromorphic
on a plane domain $D$ such that $hh^{(k)}$ has no zeros in $D$, for each
$h \in {\mathcal F}$. Then the family
$\{ h'/h : h \in {\mathcal H} \}$ is normal on $D$.
\end{thm}

\begin{lem}[\cite{Lasing2016}]
 \label{cor1}
Let $h$ be a transcendental meromorphic function in the plane such that $h'/h$ has finite lower order and $h'/h$ and $h''/h'$ have finitely many zeros.
Then $h''/h'$ is a rational function and $h$ has finite order and finitely many poles.
\end{lem}

The next two lemmas involve Tsuji's analogue \cite{Tsuji}
for the upper half-plane of Nevanlinna's characteristic function,  which was developed further by Levin and Ostrovskii \cite{LeO} 
(see also \cite{BEL,GO}). The first is directly related to Theorem \ref{FHPthm}(i) and was deduced in \cite{Laams09}  from  Frank's method \cite{FHP}.

\begin{lem}[\cite{Laams09,LaJdA2013}]
 \label{lem1}
Let $f$ be  a non-constant meromorphic function in the plane which satisfies at least one of the following two conditions:
(a)  $f$ and $f''$ have finitely many non-real zeros and poles;
(b) $f$ and $f^{(m)}$ have finitely many non-real zeros, for some $m \geq 3$. 
Then the Tsuji characteristic $T_0(r, f'/f)$ in the upper half-plane satisfies
$T_0(r, f'/f) = O( \log r )$ as $r \to \infty$.
\end{lem}
The following lemma is due to Levin and Ostrovskii \cite{LeO} (see also \cite{BEL,GO} and \cite[Lemma 2.2]{LaJdA2013}). 
\begin{lem}[\cite{LeO}]
 \label{lem2}
Let $H$ be a non-constant meromorphic function in the plane.
If $H$ and $G = \widetilde H$ satisfy, as $r \to \infty$, 
$$
N(r, H) = O( r \log r ) \quad \hbox{and} \quad 
T_0(r, H) + T_0\left(r, G \right) = O( \log r) ,
$$
then $T(r, H) = O( r \log r )$ as $r \to \infty$.
\end{lem}

\begin{lem}
 \label{lemnorfam}
Let $0 < \varepsilon < \pi/8$, $R > 0$  and  $K > 1$. 
Let $(h_n)$ be a sequence of meromorphic functions on the domain $\{ z \in \C : |z| > R, \, 0 < \arg z < \pi \}$, each of them
such that $h_n$, $h_n'$ and $h_n''$ have no zeros there. Suppose that there exists a positive sequence 
$(r_n)$ such that $r_n \to \infty$ and 
\begin{equation}
 \label{nrb0}
\min \left\{ \left| z \, \frac{h_n'(z)}{h_n(z)} \right| : K^{-1} r_n \leq |z| \leq Kr_n , \, \varepsilon \leq \arg z \leq \pi - \varepsilon \right\} \to 0
\end{equation}
as $n \to \infty$. Then 
$$
\max \left\{ \left|  z \, \frac{ h_n'(z)}{h_n(z)} \right| : K^{-1} r_n \leq |z| \leq Kr_n , \, \varepsilon \leq \arg z \leq \pi - \varepsilon \right\} \to 0
$$
as $n \to \infty$. 
\end{lem}
\textit{Proof.} 
For $q=1, 2$ let 
$$
D_q = \{  z \in \C :  K^{-q}  < |z| < K^q , \, \varepsilon /q < \arg z < \pi - \varepsilon /q \}
$$
and let $E_1$ be the closure of $D_1$. Let $n_0 \in \N$ be large.
By Theorem \ref{norfam} the 
functions $p_n(z) = r_n h_n'(r_n z)/h_n(r_nz)$, $n \geq n_0$,  form a normal family of zero-free meromorphic functions on $D_2$. Assuming that the assertion of the 
lemma is false gives, after passing to a subsequence if necessary, 
\begin{equation}
 \label{nrb1}
\liminf_{n \to \infty} \left( 
\sup \left\{ \left| p_n(z) \right| : z \in E_1  \right\} \right) >  0. 
\end{equation}
On the other hand (\ref{nrb0}) implies that there exist $u_n \in E_1$ with $\lim_{n \to \infty}  p_n( u_n   )  =  0$.
After taking a further subsequence, if necessary, it may be assumed that, as $n \to \infty$, the points  $u_n$ converge to some $ u^* \in E_1 \subseteq D_2 $ and 
the functions $p_n$ converge locally uniformly on $D_2$ to some $p$ with $p(u^*)=0$. Thus $p$ is meromorphic on $D_2$ and $p \equiv 0$ by Hurwitz' theorem, 
which contradicts (\ref{nrb1}). 
\hfill$\Box$
\vspace{.1in}

\begin{lem}
 \label{lemlogder3}
Let $ 0 < \varepsilon < \pi/8 $ and  $K> 4$, and let $q$ be a positive integer. Then there exists $C_1 > 0$ with the following property.

Let $R \geq 1$ and let
the function $h$ be meromorphic on  $D_R = \{ z \in \C : \, |z| > R, \, {\rm Im} \, z > 0 \}$, and assume that $h^{(q)}(z) \not \equiv 0$ on $D_R$,
and that
$h'$ has no zeros in $D_R$. 
Let $(r_n)$, $(\rho_n)$ and $(S_n)$
be  positive sequences such that  $\lim_{n \to \infty} r_n = \infty$ and $\lim_{n \to \infty} \rho_n = 0$. For each $n$, assume that  
\begin{equation}
 \label{nrd1}
 \max \left\{ \left| z \, \frac{h^{(q+1)}(z)}{h^{(q)}(z)} \right| : z \in \Omega_n \right\}  \leq \rho_n ,
\end{equation}
where $\Omega_n =  \left\{ z \in \C : K^{-1} r_n \leq |z| \leq Kr_n , \, \varepsilon \leq \arg z \leq \pi - \varepsilon \right\}$, and that
there exists $z_n$ with 
\begin{equation}
\label{nrd2}
\frac{ r_n }2 \leq |z_n| \leq 2r_n , \quad  2 \varepsilon \leq \arg z_n \leq \pi - 2 \varepsilon , \quad 
  \left| z_n \, \frac{ h'(z_n)}{h(z_n)} \right|  < S_n  .
\end{equation}
Then, for all sufficiently large~$n$,  the set
$$
\left\{ \theta \in  [ \varepsilon ,  \pi - \varepsilon ] : \left| t_n  \, \frac{h'(t_n e^{i \theta } )}{h(t_n e^{i \theta } )} \right| < C_1 S_n  \right\} ,
\quad t_n = K^{-1} r_n  , 
$$
has linear measure at least $\pi/2$. 
\end{lem}
\textit{Proof.} 
By (\ref{nrd1}) there exists $c = c(n) \in \C$ such that integrating from $iKr_n $ to $z \in \Omega_n$ gives   
$$
\log h^{(q)}(z) = c +  o(1) , \quad h^{(q)}(z) = e^c \left(1 +  \frac{\delta_q(z)}{q!} \right) ,  \quad \delta_q(z) = o(1).
$$
It may be assumed that $e^c = q!$, since $h'/h$ and $h^{(q+1)}/h^{(q)}$ are unchanged if $h$ is replaced on $\Omega_n$ by $q!e^{-c(n)} h$.
Thus repeated integration gives a monic polynomial $P=P_{q,n}$, of degree $q$, with the property that, for $j=0, \ldots, q$ and for all $z$ in $\Omega_n$,
$$
h^{(j)}(z) = P^{(j)}(z) + \delta_j(z), \quad \delta_j(z) = o( |z|^{q-j}), \quad \delta_j'(z) = \delta_{j+1}(z) = o( |z|^{q-j-1}); 
$$
here all these estimates hold as $n \to \infty$, uniformly on $\Omega_n$, and the last estimate for $j=q$ follows from (\ref{nrd1}). 

Let $n$ be large and denote by 
$c_j$  positive constants  which are independent of $n$. Then  (\ref{nrd2}) delivers a small  $c_0$ such that the disc $|z-z_n| < c_0 r_n$ lies in $D_R$,
and since $h' $ has no zeros in $D_R$  it follows from the minimum principle that 
\begin{equation}
 \label{nrd2a}
\min \left\{ \left| w \, \frac{ h'(w)}{h(w)} \right| : | w- z_n| = s  \right\}  \leq S_n \quad \hbox{for all} \quad  s \in [ c_0 r_n/4, c_0 r_n/2] .
\end{equation}

Now let $\{ B_j \} = \{ B_{j,n} \} $ denote the collection of all zeros of $P$ and $P'$. Let 
$Q_n$ be the closed set obtained by deleting from $\Omega_n$ the open discs $E_j$ 
of centre $B_j$ and radius $c_1 r_n$, where $c_1$ is assumed to be small. Then  $z \in Q_n$ gives  
$| z-B_j| > c_2 |z|$ 
for every $j$, and hence  
$$
|P(z)| > c_3 |z|^q , \quad |P'(z)| > c_3 |z|^{q-1} , \quad 
\left| \frac{P'(z)}{P(z)} \right|  + \left| \frac{P''(z)}{P'(z)} \right|  < \frac{c_4}{|z|} .
$$
For $z \in Q_n$ it follows that $\phi = h'/h$ satisfies 
\begin{eqnarray}
 \left| \frac{\phi'(z)}{\phi(z)} \right|  &=&  \left| \frac{P''(z)  + \delta_2(z)}{ P'(z) + \delta_1(z)} - \frac{P'(z)  + \delta_1(z)}{ P(z) + \delta_0(z)}  \right|  
\nonumber \\
&=& \left| \frac{P''(z) }{ P'(z) (1+ o(1) )} - \frac{P'(z)  }{ P(z) (1 + o(1))}  + \frac{o(1)}{|z|}  \right|  \leq  \frac{c_5}{|z|} .
\label{nrd3}
\end{eqnarray}

Provided $c_1$ was chosen small enough, the following exist:
a real number $s_n \in  [c_0 r_n/4, c_0 r_n/2] $  such that the circle $|z- z_n| = s_n$   meets none of the discs $E_j$;
a real number $u_n \in  [ r_n/3 , 3 r_n ]$ such that the circle $|z|= u_n$ meets $|z- z_n| = s_n$ but none of the  $E_j$;
a set  $T_n \subseteq [ \varepsilon ,  \pi - \varepsilon ]$, of linear measure at least $\pi/2$, such that for $\theta \in T_n$ the 
line segment given by $K^{-1} r_n \leq |z| \leq K r_n$, 
$\arg z = \theta $, lies in $Q_n$. 
Using (\ref{nrd2a}), choose $w_n$ with 
$$
|w_n -z_n| = s_n, \quad \left| w_n \, \frac{h'(w_n)}{h(w_n )}  \right|  \leq  S_n .
$$
For $v = t_n e^{i \theta } = K^{-1} r_n e^{i \theta }$ with $\theta \in T_n$, there exists a path $\Gamma_v \subseteq Q_n$, joining $v$ to $w_n$, which 
consists of part of the ray $\arg z = \theta $ and arcs of the circles $|z| = u_n$ and  $|z-z_n| = s_n$. 
The path  $\Gamma_v$  has length at most $c_6 r_n$, and so integrating $\phi'/\phi$  along $\Gamma_v$   gives, using (\ref{nrd3}), 
$$
\left| v \, \frac{h'(v)}{h(v)}  \right| = | v \phi (v)| 
 < c_7  | v \, \phi(w_n)| = c_7 \left| v \, \frac{h'(w_n)}{h(w_n)}  \right|  < c_8  \left| w_n \, \frac{h'(w_n)}{h(w_n)}  \right|  \leq c_8 S_n .
$$
\hfill$\Box$
\vspace{.1in}

\begin{lem}
 \label{lemf1}
Let $B \in \C$ with $|B|=1$ and $L \in \Z \setminus \{ -1 \}$, and let 
\begin{equation}
 \label{Fdef1}
F(z) = (e^z-1)^L (e^z-B).
\end{equation}
If all zeros of $F''$ lie on $i\R = \{ ix : x \in \R \}$ then $F$ is given by one of the following:\\
$$
(i) \quad F(z) = e^{2z} - 1; \quad (ii) \quad F(z) = e^z-B; \quad (iii) \quad F(z) = \frac1{e^z-1} . 
$$
\end{lem}
\textit{Proof.} Note first that $L=0$ gives (ii) immediately, while if $L=1$ then $F''(z) = 4e^{2z} - (1+B)e^z$, which has zeros off $i \R$ unless
$B=-1$, in which case $F$ is given by (i). Assume henceforth that $L \neq 0, \pm 1$, so that $LB + 1 \neq 0$, since $|B|=1$. Now  write $X = e^z$ and 
\begin{eqnarray*}
 F(z) &=& P(X) = (X-1)^L (X-B),\\
 F''(z) &=& XP'(X) + X^2 P''(X) \\
&=& X (L (X-1)^{L-1} (X-B) +  (X-1)^L ) + \\
& & +  X^2 (  L(L-1) (X-1)^{L-2} (X-B) +  2L (X-1)^{L-1} ),
\end{eqnarray*}
from which it follows that 
\begin{eqnarray*}
Q(X) &=& 
 X^{-1}(X-1)^{2-L} F''(z) \\
&=& X^2 (L+1)^2 + X ( - 3L -2 -L^2 B ) + LB+1 
\end{eqnarray*}
satisfies $Q(0) = LB + 1 \neq 0$. 
If $Q(C) =0$ and $C \neq 1$ then $e^z = C$ implies $F''(z)=0$. 
Hence the fact that all zeros of $F''$ lie on $i \R$ forces each root of $Q$ to have modulus $1$, so that 
$$|L+1|^2 = |LB+1| \leq |L| + 1,$$
which is impossible if $L \geq 2$.  
Now suppose that $L=-n \leq -2$. Then $(n-1)^2 \leq n+1$ and so $n \leq 3$,
giving $L = -2, -3$. Now $L=-2$ forces $|-2B+1| = 1$ and so $B=1$, which leads to (iii).
Finally, if $L=-3$ then $|-3B+1| = 4$, from which it follows that $B=-1$ and 
$$
Q(X) = 4X^2 + 16X + 4 ,
$$
which does not have two roots of modulus $1$. 
\hfill$\Box$
\vspace{.1in}

\begin{lem}
\label{lemtan}
 Let $A \in \C$. Then all solutions of $ \tan z = z - A $ are real if and only if $A \in \pi \Z = \{ n \pi: n \in \Z\}$.
\end{lem}
\textit{Proof.} This is proved in \cite{HinRos} (starting from formula (5) on p.73, with $B = 1$ in the notation there). 
However, since the method  in  \cite{HinRos} uses \cite[Lemma 8]{HelW1}, the proof of which is lengthy, it seems worth including the following
self-contained argument. 
First, suppose that $p(z) = \tan z$ has a fixpoint $w$ in the open upper half-plane $H^+$. Since $p$ maps $H^+$ into itself, but not univalently, Schwarz' lemma
implies that $|p'(w)| < 1$ and that the iterates of $p$ converge to $w$ on $H^+$. 
Hence $w$
must lie on the positive imaginary axis $i \R^+$, since $p( i \R^+) \subseteq i \R^+$, which contradicts the fact
that simple estimates give $\tanh y < y$ for $y > 0$. Hence all fixpoints of $\tan z$ are real and periodicity implies that so are all solutions 
of  $ \tan z = z - A $ for  $A \in \pi \Z$. 


Now suppose that all solutions of $ \tan z = z - A $ are real. 
The  real meromorphic function $g(z) = z - \tan z$ has no finite asymptotic values, and so  no Picard values: thus $A  \in \R$.
Suppose that $n \in \Z$ and $n \pi < A < (n+1) \pi$. Then $g(x)$ is decreasing on the interval 
$((n+1/2) \pi , (n+3/2) \pi ) $ and has a fixpoint at $(n+1) \pi$, which is a zero of $g - (n+1) \pi$ of multiplicity $3$. 
Hence there exists a level curve $\gamma$ on which $g(z)$ is real and decreasing, which starts at $(n+1) \pi$ and enters the upper half-plane. 
Since $g$ has no finite asymptotic values, and all critical points of $g$  are  fixpoints of $g$ in $\pi \Z$, the curve $\gamma$ must pass through a non-real 
$A$-point of $g$. 
\hfill$\Box$
\vspace{.1in}

\section{Intermediate steps for Theorems \ref{thm1} and \ref{thm2}}\label{inter}

Throughout this section, let $M \geq 4$ be an integer and let $f$ be  a strictly non-real meromorphic function  in the plane such  that   
$f, f', \ldots, f^{(M+3)}$ all have finitely many zeros in $\C \setminus \R$. 

\begin{lem}
 \label{lemfpoly}
If $f$ is a polynomial then  $f'/f$ is  rational and 
$f$ has at least one non-real zero. 
\end{lem}
\textit{Proof.} This follows at once from $f$ being strictly non-real. 
\hfill$\Box$
\vspace{.1in}

Assume henceforth that $f$ is not a polynomial, and let $g(z) = \widetilde f (z) = \overline{ f( \bar z)} $. Then
Lemma~\ref{lem1} shows that the Tsuji characteristics \cite{BEL,GO,Tsuji} of $f'/f$ and $g'/g$ satisfy 
\begin{equation}
 \label{firstest}
T_0(r, f'/f) + T_0(r, g'/g) = O( \log r ) \quad \hbox{as} \quad r \to \infty . 
\end{equation}
The lemma of the logarithmic derivative for the Tsuji characteristic \cite{LeO} and the formulas 
\begin{equation}
 \label{well}
\psi = \frac{\phi'}{\phi} , \quad 
\frac{\phi''}{\phi'} = \frac{\psi'}{\psi} + \psi ,
\end{equation}
then deliver, for  all $m \geq 0$,  
\begin{equation}
 \label{nrc2}
T_0(r, f^{(m+1)}/f^{(m)}) + T_0(r, g^{(m+1)}/g^{(m)}) = O( \log r )  \quad \hbox{as} \quad r \to \infty .
\end{equation}
For $0 \leq m \leq M+1$ write 
\begin{equation}
 \label{nrc1}
F_m(z) = z - \frac{f^{(m)}(z)}{f^{(m+1)}(z)}, \quad G_m(z) =  z - \frac{g^{(m)}(z)}{g^{(m+1)}(z)} = \widetilde F_m(z) .
\end{equation}


\begin{lem}
 \label{Kmexist}
Let $0 \leq m \leq M+1$.
Then the functions $F_m$ and $G_m$ are non-constant, and
there exists a  meromorphic function $K_m$, with finitely many zeros and poles, such that 
\begin{equation}
 \label{nrc3}
 F_m' = \frac{f^{(m)}f^{(m+2)}}{(f^{(m+1)})^2 } = K_m \left( \frac{g^{(m)}g^{(m+2)}}{(g^{(m+1)})^2 } \right) = K_m G_m'   .   
\end{equation}
The function $K_m$  satisfies $|K_m(x)| = 1$ for all $x \in \R$ and 
there exist a rational function $R_m$ and a real number $a_m$ 
such that 
\begin{equation}
\label{nra1}
K_m(z) = R_m(z) e^{i a_m z} .
\end{equation}
Furthermore, if $f^{(m)}$, $f^{(m+1)}$ and $f^{(m+2)}$ have only real zeros, then $R_m$ is constant. 
\end{lem}
\textit{Proof.} The first assertion holds since if $F_m$ is constant then $F_m' $ and $f^{(m+2)}$ vanish identically.
Now $K_m$ has finitely many zeros and poles, since $f, \ldots, f^{(M+3)} $ have finitely many non-real zeros,
and $\widetilde K_m = 1/K_m$. Finally, 
(\ref{nrc2}) and Lemma \ref{lem2} imply that 
(\ref{nra1}) holds.
\hfill$\Box$
\vspace{.1in}

\begin{lem}
 \label{lemnosimple}
For  $0 \leq m \leq M+1$: \\
(a) every real multiple zero of $f^{(m)}$ is a $1$-point of $K_m$; \\
(b) if $K_m$ is constant in (\ref{nrc3}), then either $F_m = G_m$ or $f^{(m)}$ has at most one real zero, counting multiplicities;\\
(c) every real simple zero $a$ of $f^{(m+1)}$ either is a multiple zero of $f^{(m)}$ or satisfies $K_m'(a)=0$.
\end{lem}
\textit{Proof.} To prove (a) and (b) take a real zero $x_0$ of $f^{(m)}$ of multiplicity $p$. Then $x_0$ is a zero of $g^{(m)}$ of the same multiplicity,
and a common fixpoint of $F_m$ and $ G_m $.
If $p \geq 2$ then 
$$
F_m'(x_0) = G_m'(x_0) = \frac{p-1}{p} , \quad K_m(x_0) = 1,
$$
which proves (a). Next, if 
$K_m$ is constant but $F_m \neq G_m$  then there exists   $c_m \in \C$ such that 
$$
  G_m \neq F_m = K_m G_m + c_m  ,
$$
so that  $F_m$ and $ G_m $ have at most one common fixpoint, and none at all if $K_m=1$. In view of (a), this proves  (b). 

To prove (c) take a real simple zero $a$ of $f^{(m+1)}$ which is not a zero of $f^{(m)}$. Then  $a$ is a simple pole of $F_m$, 
and there exists $b \in \C \setminus \{ 0 \} $ such that, as $z \to a$, 
$$
F_m(z) = \frac{b}{z-a} + O(1), \quad  F_m'(z) = \frac{-b}{(z-a)^2} + O(1) , \quad  G_m'(z) = \frac{- \bar b}{(z-a)^2} + O(1) .
$$
This implies that $K_m'(a) =0$. 
\hfill$\Box$
\vspace{.1in}

The next three lemmas will treat  a number of special cases. 

\begin{lem}
 \label{lem4}
Assume that  $0 \leq m \leq M+1$ and at least one of the following  holds:\\
(i) $f^{(m+1)}/f^{(m)}$ is real meromorphic;\\
(ii) $F_m = G_m$;\\
(iii) $g^{(m)} = c_m f^{(m)}$ for some  $c_m \in \C $. \\
Then  $f$ is a rational
function with at least one non-real zero.  
\end{lem}
\textit{Proof.} It is clear from (\ref{nrc1}) that (i) implies (ii) and (ii) implies (iii). Assume therefore that (iii) holds:
then $|c_m | = 1$, because  $g = \widetilde f$ and $f^{(m)} \not \equiv 0$. 
Moreover, $m \geq 1$, since $f$ is strictly non-real,
and $f$ and $g$ have the same poles with the same multiplicities.  Hence there exists a non-constant meromorphic function $H$ 
with finitely many zeros and poles such that,
using (\ref{firstest}), 
\begin{equation}
 \label{nrc1a}
 g = H f, 
\quad \widetilde H = \frac1{H}, 
\quad \frac{g'}{g} - \frac{f'}{f} = h = \frac{H'}{H} , \quad T_0(r,h)  = O( \log r )  \quad \hbox{as} \quad r \to \infty .
\end{equation}
Furthermore, integration gives 
a polynomial $P \not \equiv 0$, of degree at most $m - 1 \leq M$, with 
\begin{equation}
 \label{nrc4a}
g = H f = P + c_m f, \quad f = \frac{P}{H - c_m} , \quad \frac{f'}{f} = \frac{P'}{P} - \frac{H'}{H-c_m} = \frac{P'}{P} - \frac{h }{1 - c_m H^{-1}}.
\end{equation}
Hence $T_0(r, H ) = O( \log r)$ and $T(r, H) = O(r \log r )$ as $r \to \infty$, by (\ref{firstest}), (\ref{nrc1a}) and  Lemma~\ref{lem2}.
Thus $H(z) = T_1(z) e^{i a_1 z } $, where $a_1 \in \R$ and $T_1$ is a rational function with $|T_1(x)| = 1$ on $\R$. 

If $H$ is transcendental then (\ref{nrc1a}) and (\ref{nrc4a}) show that  $f$ satisfies the hypotheses 
of \cite[Lemma 2.5]{LaJdA2013}, and so $f'''$ has infinitely many non-real zeros,
contrary to assumption. Therefore  $H$ is a rational function and so is $f$. Because  $H$ is non-constant
and $\widetilde H = 1/H$, the function  $H$ has at least one pole and, since 
$f$ and $g$ have the same poles, $f$ has at least one non-real zero. 
\hfill$\Box$
\vspace{.1in}

\begin{lem}
 \label{lemfo}
Assume that $f$ has finite order and finitely many poles. Then either $f'/f$ is a rational function, or 
$f$ satisfies (\ref{thm1con0}).
\end{lem}
\textit{Proof.} 
The hypotheses imply that there exist meromorphic functions $H$ and $K$, each with finitely many zeros and poles,  such that 
\begin{equation}
 \label{fo1}
g = Hf, \quad g' = HK f', \quad \widetilde H = \frac1{H}, \quad \widetilde K = \frac1K .
\end{equation}
Since $f$ is strictly non-real, $H$ is non-constant. 
Write
\begin{equation}
 \label{fo3}
h = \frac{H'}{H}, \quad k = \frac{K'}{K}, \quad \widetilde h = -h, \quad \widetilde k = - k.
\end{equation}
Then $h$ and $k$ have finitely many poles and so are rational functions, since $f$ has finite order.
Moreover, $h$ does not vanish identically, since $H$ is non-constant, and $h'/h$ is real. 
 
Now (\ref{fo1}) and (\ref{fo3}) yield
$$
g' = hHf + Hf' = HKf'.
$$
Here $K-1$ cannot vanish identically because $h$ does not. It follows that 
\begin{equation}
 \label{fo4}
L = \frac{f'}{f} = \frac{h}{K-1} .
\end{equation}
If $K$ is a rational function, then so is $f'/f$.

Assume henceforth that  $K$, which  has finitely many zeros and poles,
is transcendental; then $k \not \equiv 0$ in (\ref{fo3}). Moreover, Lemma \ref{lem2}, (\ref{firstest}), (\ref{fo1}) and (\ref{fo4}) imply 
that $K(z) = T_1(z)e^{i A_1z} $, where $A_1 \in \R \setminus \{ 0 \}$ and $T_1 $ is rational with $|T_1(x)| = 1$ on $\R$. 

If $h = \pm k$ then (\ref{fo4}) shows that 
$$
 \frac{f'}{f} = \pm \frac{K'}{K(K-1)} 
, \quad f = c (1-1/K)^{\pm 1}, \quad c \in \C \setminus \{ 0 \},
$$
and so $f$, which has finitely many poles,
must satisfy (\ref{thm1con0}). 

Assume henceforth 
that $h \neq \pm k$. 
Combining (\ref{well}), (\ref{fo3}) and (\ref{fo4}) leads to 
\begin{equation}
 \label{fo5}
\frac{f''}{f'} = L + \frac{L'}{L} = \frac{h}{K-1} + \frac{h'}{h} - \frac{kK}{K-1} 
= \frac{h-k}{K-1} + \frac{h'}{h} - k  .
\end{equation}
Observe that
none of the functions 
$k \pm h'/h $, $h \pm h'/h$  vanishes identically, since $h'/h$ is real but $h$ and $k$ are not.
If $|z|$ is large, then (\ref{fo5}) shows that $z$ is a zero of $f''/f'$ if and only if $z$ is a solution of the following equations:
$$
\frac{h-k}{K-1} = k - \frac{h'}{h} ; \quad K-1 = \frac{h-k}{k-h'/h} ; \quad K = \frac{h - h'/h}{k-h'/h}.
$$
Thus   $f''/f'$ has 
infinitely many real zeros $x$ which satisfy, by (\ref{fo1}) and (\ref{fo3}), 
$$
  \frac{k(x)-h'(x)/h(x)}{h(x) - h'(x)/h(x)} = \frac1{K(x)} = \overline{ K(x) } = \frac{-h(x) - h'(x)/h(x)}{-k(x)-h'(x)/h(x)} = \frac{h(x) + h'(x)/h(x)}{k(x)+h'(x)/h(x)}.
$$
Because $k$ and $h$ are rational functions, this forces 
$$
k^2 - (h'/h)^2  = h^2 - (h'/h)^2 , \quad  h^2 = k^2,
$$ 
contradicting the assumption that $h \neq \pm k$. 
\hfill$\Box$
\vspace{.1in}

\begin{lem}
 \label{lemf'frational3}
Assume that either $f'/f$ is a rational function or 
$f$ satisfies (\ref{thm1con0}), and that $f$, $f'$, $f''$ and $ f'''$ have only real zeros.  
Then $f$ is given by (\ref{thm1con1}).
\end{lem}
\textit{Proof.} Suppose first that $f$ satisfies (\ref{thm1con0}).
Then $f''/f'$ is a rational function, and so is $F_1$ in (\ref{nrc1}). 
Moreover,  the function
$K_1$  in (\ref{nrc3}) is rational and free of  zeros and  poles, and so is constant, but  Lemma \ref{lem4} implies that $G_1 \neq F_1$.
Applying Lemma \ref{lemnosimple} shows that $f'$ has at most one zero, and that any  zero of $f'$ is real and simple. Now (\ref{thm1con0}) gives
$$
\frac{f'}{f-B} = \frac{T'}{T} + iA , \quad A \neq 0.
$$
If $T$ is non-constant then $T'/T$ has at least two poles in $\C$, since $\widetilde T = 1/T$, and so $f'$ has at least two zeros in $\C$, counting multiplicities.
This is a contradiction, and so 
$f$ is given by (\ref{thm1con1})(iii). 

Assume henceforth that $R = f'/f$ is a rational function. Then so are $F_0$ and $F_1$  in (\ref{nrc1}),
and the same argument as in the previous paragraph shows that $K_0$ and $K_1$ are constant.
However,   Lemmas \ref{lemnosimple} and \ref{lem4} and 
the fact that $f$ is strictly non-real imply the following: $G_0 \neq F_0$ and $G_1 \neq F_1$; neither $f'/f$ nor $f''/f'$ is real;
any  zero of $f$ is real, simple and unique, and the same  applies to zeros of $f'$ and~$R$. 

Suppose first that  $R(\infty) = \infty$. Then $R$,  since it has at most one zero, 
must have form $R(z) = \alpha (z - x_0)$ with $x_0 \in \R$ and $0 \neq \alpha \in \C$, so that 
$f''(z)/f(z) = \alpha + \alpha^2 (z-x_0)^2$. Because $f''$ has only real zeros, $\alpha $ is real and so is $f'/f$, a contradiction. 

If $R$ is  a non-zero constant, then $f$ satisfies~(\ref{thm1con1})(iii). Suppose next that $R$ is non-constant, with $R(\infty) \neq 0, \infty$. 
Then $R$ is a M\"obius transformation, since it has at most one zero. Applying a change of variables $w = a_1 z + b_1$ with $a_1, b_1 \in \R$
makes it possible to assume that the unique zero of $R$ is at the origin, and that
$$
\frac{f'(z)}{f(z)} = 
R(z) = \frac{az}{z - z_0} = a + \frac{az_0}{z-z_0} , \quad 
\frac{f''(z)}{f(z)}  = \frac{a^2 z^2 - az_0}{(z-z_0)^2} ,
$$
where $a, z_0 \in \C \setminus \{ 0 \}$.
Here $b = az_0 $ is an integer and $z_0 \not \in \R$, since otherwise $a$ and $f'/f$ are real. Thus $b$ must be negative and $z_0$ is a pole of $f$ 
and a double pole of $f''/f$. 
Next, 
$z_0/a $ must be real and positive, since $f''$ has only real zeros, and so must $-z_0^2 = -bz_0/a$.
Now write 
$$
\frac{f'''(z)}{f(z)} =  \frac{az(a^2 z^2 - az_0)}{(z-z_0)^3} + \frac{2a^2 z}{(z-z_0)^2} - \frac{2(a^2 z^2 - az_0)}{(z-z_0)^3}  =
\frac{a( a^2 z^3  - 3az_0 z + 2z_0)}{(z-z_0)^3} .
$$
Because $z_0^2$
and $az_0$ are real, so is $a^2$. Since $z_0$ is not real, $f'''$ must have at least one non-real zero, a contradiction.

Assume next that $R$ has a simple zero at infinity. If $R$ has no zeros in $\C$ then $f/f' = 1/R$ is a linear polynomial and $f$ satisfies~(\ref{thm1con1})(i).
If $R$ has a zero in $\C$ then it has exactly one zero and two poles there, and it may be assumed that 
\begin{equation}
 \label{Rforma}
\frac{f'(z)}{f(z)} = 
R(z) = \frac{az}{(z-z_1)(z-z_2)} , \quad a, z_1, z_2 \in \C \setminus \{ 0 \}, \quad z_1 \neq z_2.
\end{equation}
Here the residues $r_1=az_1/(z_1-z_2)$ and $ r_2= a z_2/(z_2-z_1)$ must be integers, and $r_1/r_2 = - z_1/z_2 $ is real. If either residue $r_j$ is positive, then
$z_1$ or $z_2$ is real, so that both are real, and so is $a$, contradicting the fact that $R = f'/f$ is not real. So both $r_j$ are real and negative,
as are $z_1/z_2$ and  $a$, and $f(z_1) = f(z_2) = \infty$. Now
\begin{eqnarray*}
 \frac{f''(z)}{f(z)} &=& \frac{a^2z^2}{(z-z_1)^2(z-z_2)^2} + \frac{a}{(z-z_1)(z-z_2)} - \frac{az}{(z-z_1)^2(z-z_2)} - \frac{az}{(z-z_1)(z-z_2)^2} \\
&=&  \frac{a^2 z^2 + a (z-z_1)(z-z_2 ) - az (z-z_2) - az (z-z_1)} { (z-z_1)^2(z-z_2)^2 } = \frac{(a^2-a) z^2 + az_1 z_2 }{(z-z_1)^2(z-z_2)^2 } .
\end{eqnarray*}
Since $a < 0$ this forces $z_1 z_2 $ to be real and positive, and so $z_1^2 $ and $z_2^2$ are real and negative. Next,
\begin{eqnarray*}
 \frac{f'''(z)}{f(z)} &=&  \frac{az((a^2-a) z^2 + az_1 z_2 )}{(z-z_1)^3(z-z_2)^3 } + \frac{(a^2-a) 2z  }{(z-z_1)^2(z-z_2)^2 } + \\
& & - \frac{2((a^2-a) z^2 + az_1 z_2 ) }{(z-z_1)^3(z-z_2)^2 }  - \frac{2((a^2-a) z^2 + az_1 z_2 ) }{(z-z_1)^2(z-z_2)^3 } \\
&=& \frac{((a^2-a) z^2 + az_1 z_2 )(az -4z + 2(z_1+z_2)) + (a^2-a) 2z (z-z_1)(z-z_2) }{(z-z_1)^3(z-z_2)^3 } \\
&=& \frac{a(a-1)(a-2) z^3 + z_1 z_2 (3 a^2 - 6a) z + 2a z_1 z_2 (z_1 + z_2) }{(z-z_1)^3(z-z_2)^3 } .
\end{eqnarray*}
But $a < 0$, and $f'''/f$ has triple poles at $z_1$ and $z_2$. 
Hence $f'''/f$ has three zeros in $\C$, counting multiplicities, all of them real. Because $z_1z_2$ is real, $z_1+z_2$ must be real, and so $0$.
But then (\ref{Rforma}) implies that $f'/f$ is real, a contradiction.

Finally, suppose that  $R$ has  a zero at $\infty$ of multiplicity at least two. Then integration of $R$ around a circle $|z| = r $ with $r$ large  
shows that $f$ has  in $\C$ the same number of zeros as poles, counting multiplicities, and so exactly one of each. Hence  $f$ 
satisfies (\ref{thm1con1})(ii).

\hfill$\Box$
\vspace{.1in}

Assume for the remainder of this section that $f$ has  either infinite order of growth or infinitely many poles. 
Then $f^{(m+1)}/f^{(m)}$ is transcendental,  for each $m \geq 0$.

\begin{lem}
 \label{lemFmFm+1}
The following statements all hold.\\
(i) If $0 \leq m \leq M+1$ and $K_m$ is  constant then $f^{(m)}$ has finitely many zeros.\\
(ii) If $0 \leq m \leq M$ and  $K_m$ and $K_{m+1}$ are both non-constant, then $\overline{N}(r, 1/f^{(m+1)}) = O(r)$ as $r \to \infty$. \\
(iii) If $0 \leq m \leq M$ and  $K_m$ and $K_{m+1}$ are both non-constant rational functions, 
then $f^{(m+1)}$ has finitely many zeros;\\
(iv) If $0 \leq m \leq M$ and  $K_m$ and $K_{m+1}$ are both  rational functions, 
then $f^{(m)}$ or $f^{(m+1)}$ has finitely many zeros.
\end{lem}
\textit{Proof.} Since $f$ is transcendental by assumption, Lemma \ref{lem4} shows that $F_m \neq G_m$ for $0 \leq m \leq M+1$. 
Thus (i) follows from Lemma \ref{lemnosimple}(b).

Next, assume the hypotheses of (ii),
and let $x_0$ be a real zero of $f^{(m+1)}$.
By Lemma \ref{lemnosimple},  either $x_0$  is  a multiple zero of $f^{(m)}$ or 
$f^{(m+1)}$, and hence  a $1$-point of
$K_m$ or $K_{m+1}$, or $x_0$ is a zero of $K_m'$. Now  (ii) and (iii) follow, by (\ref{nra1}), and combining (i) and (iii) gives (iv). 
 \hfill$\Box$
\vspace{.1in}

\begin{lem}
 \label{lemM=1}
There exists $\alpha > 0$ such that, for $1 \leq m \leq M+2$,
\begin{equation}
 \label{nrc2b}
T(r, f^{(m+1)}/f^{(m)}) +  T(r, g^{(m+1)}/g^{(m)}) < \alpha r  \quad \hbox{as} \quad r \to \infty .
\end{equation}
\end{lem}
\textit{Proof.} 
If $K_0$ or $K_1$ is constant, then $f$ or $f'$ has finitely many zeros, by Lemma \ref{lemFmFm+1}.
If  $K_0$ and $K_1$ are both non-constant then $\overline{N}(r, 1/f') = O(r)$ as $r \to \infty$. This implies that 
\begin{equation}
 \label{zerosest}
\overline{N}(r, 1/f^{(m)}) = O(r) \quad \hbox{as} \quad r \to \infty 
\end{equation}
holds for $m= 0$ or $m=1$. Since $M \geq 4$, the same argument may be applied to $K_4$ and $K_5$ to show that  (\ref{zerosest}) holds 
for $m= 4$ or $m=5$. 
This delivers $p \in \{0, 1\}$ and $q \in \{ 3, 4, 5 \}$ such that (\ref{zerosest}) holds for $m=p$ and $m=p+q$. 
Now Theorem \ref{FHPthm}  implies that there exists $d_1 > 0$ with 
$$
T(r, f^{(p+1)}/f^{(p)}) \leq d_1 \left( \overline{N}(r, 1/f^{(p)}) +  \overline{N}(r, 1/f^{(p+q)}) \right)  + O( \log r ) = O(r) 
$$
as $r \to \infty$ outside a set of finite measure. This gives (\ref{nrc2b}) for some $m \in \{ 0, 1 \}$ and positive~$\alpha$. The existence of
$\alpha > 0$ such that (\ref{nrc2b}) holds for $1 \leq m \leq M+2$
then follows from (\ref{well}).
\hfill$\Box$
\vspace{.1in}

\begin{lem}
 \label{lemonenozeros}
Call an integer $m $ exceptional if $m \geq 0$ and $f^{(m+1)}$ has finitely many zeros. 
Then there exists at most one exceptional $m$. 
\end{lem}
\textit{Proof.} 
Suppose that there exist $m $
and $m'$ with $0 \leq m < m' $ such that $f^{(m+1)}$ and
$f^{(m'+1)}$ have finitely many zeros.
If $m' \geq m+2$ then $f^{(m+1)}$ has finite order and finitely many poles,
by Theorem \ref{FHPthm}, a contradiction. 
If $m'=m+1$ then the same contradiction is obtained by applying Lemma  \ref{cor1} 
to $f^{(m)}$, 
using the fact that $f^{(m+1)}/f^{(m)} $ has finitely many zeros and hence 
finite order by (\ref{nrc2}) and  Lemma \ref{lem2}.
\hfill$\Box$
\vspace{.1in}

\begin{lem}
 \label{lemexceptional}
Assume that  $m \in \{ 1, \ldots, M-2 \}$ is exceptional. 
Then
 \begin{equation}
 \label{noto(r)}
 \alpha^*  =
 \liminf_{r \to \infty} \frac{T(r, f^{(m+1)}/f^{(m)})}{r}   \in (0, + \infty) . 
\end{equation}
\end{lem}
\textit{Proof.} 
Assume  that $m$  is as in the statement  but  (\ref{noto(r)}) fails.
In view of (\ref{nrc2b}) it must therefore be the case that $\alpha^*  = 0$ in
(\ref{noto(r)}). Then there exists  a
sequence $s_n \to \infty$ with $T(s_n, f^{(m+1)}/f^{(m)}) = o(s_n)$ and so,
by (\ref{well}), (\ref{nrc3}), 
(\ref{nrc2b}) and the lemma of the logarithmic derivative,
$$
\sum_{j=m}^{m+4} T(s_n, f^{(j+1)}/f^{(j)}) = o(s_n), \quad 
T(s_n, F_{m+2}') +  T(s_n, F_{m+3}') = o(s_n).
$$
Since $m+2 \leq M$, it follows that $a_{m+2} = a_{m+3} = 0$ in (\ref{nra1}) and hence, by  Lemma \ref{lemFmFm+1},
that 
$f^{(m+2)}$ or $f^{(m+3)}$ has finitely many zeros, which contradicts
Lemma \ref{lemonenozeros}.
\hfill$\Box$
\vspace{.1in}

\begin{lem}
 \label{lemsequence}
There exist  a  real number $M_1 > 1$ and an increasing positive
sequence $(r_n )$  with limit $\infty$ such that, for all large $n$ and all $m \in \{ 1, \ldots, M \}$, 
$$
T(2r_n, f^{(m+1)}/f^{(m)} ) \leq M_1  \, T(r_n, f^{(m+1)}/f^{(m)} ) .
$$
\end{lem}
\textit{Proof.} Let $M_2 > 1$. By (\ref{nrc2b}) and a growth lemma of Hayman \cite[Lemma 4]{Hay3}, each set 
$$
X_m = \{ r \geq 1 : \, T(2r, f^{(m+1)}/f^{(m)} )> M_2  \, T(r, f^{(m+1)}/f^{(m)} ) \}
$$
has upper logarithmic density at most $d_0 = \displaystyle{\frac{\log 2}{\log M_2 } }$. Hence it suffices to take $M_2$ so large
that $Md_0 < 1$, and choose a sequence $r_n \to \infty$ in the complement of the union of the $X_m$. 
\hfill$\Box$
\vspace{.1in}

\begin{lem}
 \label{lemseq}
Let $m \in \{ 1, \ldots, M \}$ and let $\varepsilon > 0$. \\
(A) If $m$  is exceptional and  $\varepsilon $ is  sufficiently small,
then for each large $n$  there exist $L_{m,n} \in  \{ f^{(m)},  g^{(m)} \}$ and $v_n$ satisfying
\begin{equation}
 \label{nra1b1}
|v_n| = r_n, \quad 2 \varepsilon \leq \arg v_n \leq \pi - 2 \varepsilon, \quad 
\left| v_n \, \frac{L_{m,n}'(v_n)}{L_{m,n}(v_n)} \right|  < \exp(  - (1/4)  T(r_n,f^{(m+1)}/f^{(m)}) ) .
\end{equation}
(B) 
If $a_m \neq 0$ in (\ref{nra1}) then for each  large $n$ 
there exist $L_{m,n} \in  \{ f^{(m)}, f^{(m+1)} , g^{(m)}, g^{(m+1)} \}$ and $v_n$ satisfying
\begin{equation}
 \label{nra1b}
 |v_n - ir_n| < 4, \quad 
\left| v_n \, \frac{L_{m,n}'(v_n)}{L_{m,n}(v_n)} \right| < e^{ - |a_m|  r_n /8 } .
\end{equation}
\end{lem}
\textit{Proof.} 
To prove (A) assume that $m$ is exceptional
and let $n$ be large. 
Since $ f^{(m+1)}/f^{(m)}$ is transcendental and has finitely many zeros, combining Lemma \ref{lemsequence} with  a well known estimate
due to Edrei and Fuchs 
\cite[p.322]{EF2} shows that, provided $\varepsilon$ is small enough, 
the set 
$$\left\{ \theta \in [0, 2 \pi] : \, 
\left| \frac{f^{(m+1)}(r_n e^{i \theta })}{  f^{(m)} (r_n e^{i \theta }) } \right| < 
\exp( - (1/2) T(r_n, f^{(m+1)}/f^{(m)}) ) \right\}$$
has  measure at least
$16 \varepsilon$.  Hence there exist $L_{m,n} \in  \{ f^{(m)},  g^{(m)} \}$ and $v_n$
such that (\ref{nra1b1}) holds.

To prove (B), assume that $a_m \neq 0$  in (\ref{nra1}) and again
let $n$ be large. By interchanging $f$ and $g$ it may be assumed that $a_m > 0$. This implies that 
\begin{equation}
 \label{nra2}
\left| K_m(z)  \right| < e^{- |a_m| r_n/2 }   \quad \hbox{for $|z-ir_n| < 2$.}
\end{equation}
It follows immediately from (\ref{nrc3}) and (\ref{nra2}) that, for each $z$ with $|z-ir_n| < 2$,  
\begin{equation}
 \label{nra3}
\hbox{either} \quad (a) \quad |G_m'(z) | > e^{ |a_m| r_n/4 } \quad \hbox{or} \quad (b) \quad |F_m'(z) | < e^{- |a_m| r_n/4 } .
\end{equation}
Suppose first that  case (a) holds in (\ref{nra3}), for some $z$ with $|z-ir_n| < 2$. Because $n$ is large and $G_m$ has finitely many non-real poles, by  (\ref{nrc1}),
Cauchy's estimate for derivatives  implies that there exists $v_n$ with $|v_n - ir_n | <  4$ such that 
\begin{equation*}
 \label{vnest1}
|G_m(v_n)| >  e^{ |a_m| r_n/6 } , \quad \left| v_n \frac{g^{(m+1)}(v_n)}{g^{(m)}(v_n)} \right|   <     e^{ - |a_m| r_n/8 } .
\end{equation*}
This gives (\ref{nra1b}) with $L_{m,n} = g^{(m)}$.

Now suppose that  case (b) holds in (\ref{nra3}), for all $z$ with $|z-ir_n| < 2$, in which case 
\begin{equation}
 \label{Fm'est}
|F_m(z)- F_m(ir_n) | < 2 e^{ - |a_m| r_n/4 }  \quad \hbox{for $|z-ir_n| < 2$.}
\end{equation}
Choose $v_n$  with $|v_n -ir_n| < 2$  such that $| v_n - F(ir_n)| \geq 1$. Then (\ref{nrc1}) and (\ref{Fm'est}) imply that
$$
\left| \frac{f^{(m)}(v_n)}{f^{(m+1)}(v_n)} \right| = | (v_n - F_m(ir_n)) - (F_m(v_n)-F_m(ir_n)) |  \geq  \frac12 .
$$
It now follows from (\ref{nra3}) that 
\begin{equation*}
 \label{vnest2}
\left| v_n \frac{f^{(m+2)}(v_n)}{f^{(m+1)}(v_n)} \right| = \left| v_n F_m'(v_n) \cdot  \frac{f^{(m+1)}(v_n)}{f^{(m)}(v_n)} \right|
<  e^{ - |a_m| r_n/8 } .
\end{equation*}
Thus (\ref{nra1b}) holds with $L_{m,n} = f^{(m+1)}$. 
\hfill$\Box$
\vspace{.1in}

\begin{lem}
 \label{lemnomj}
Let $Y$ be the set of integers $m \in  \{ 1, \ldots, M \}$ such that 
either $m$ is exceptional or 
$a_{m} \neq 0$ in (\ref{nra1}).
Then there do not exist 
integers $m_1, m_2, m_3 \in Y$ satisfying 
\begin{equation}
 \label{mjchoose}
m_2 \geq m_1 + 2, \quad m_3 \geq m_2 + 2.
\end{equation}
\end{lem}
\textit{Proof.} Assume that  $m_1, m_2, m_3 \in Y$ satisfy (\ref{mjchoose}). 
If any $m_\nu$ is exceptional and  $m_\nu \leq M-2$ then it is unique, by Lemma \ref{lemonenozeros}: in this case 
let $\alpha^*  $ be as in (\ref{noto(r)}), so that $\alpha^* > 0$. On the other hand, if no such $m_\nu$ exists
let $\alpha^*  = 0$. In either case the set 
$\{ \alpha^* , |a_{m_1}|,  |a_{m_2}|, |a_{m_3}| \}$ has a least positive member, which will be denoted by $\beta$. 

Define $S_n$ and $K$ by
\begin{equation}
 \label{choosek}
 S_n =   e^{ - \beta  r_n /8 }, \quad 
K = 4 + \frac{128 \alpha}{\beta}   ,
\end{equation}
where $\alpha$ is as in (\ref{nrc2b}), and let $\varepsilon$ be small and positive. 

Apply  Lemma \ref{lemseq} with $m=m_\nu$ and $\nu = 1, 2, 3$;
this is possible since if $a_{m_\nu} = 0$ then $m_\nu$ is exceptional, by the definition of $Y$.
Passing to a subsequence then allows the following  assumptions to be made for $\nu= 1, 2, 3$ and all sufficiently large $n$:
first,  if $m_\nu$ is exceptional
then  (\ref{nra1b1}) holds for $m=m_\nu$, while   if  $m_\nu$ is not exceptional then $a_{m_\nu} \neq 0$ and (\ref{nra1b}) holds for $m=m_\nu$; 
second, 
$H_\nu = L_{m_\nu,n}$ is for each $n$ 
the same element of the set $\{ f^{(m_\nu)}, f^{(m_\nu+1)} , g^{(m_\nu)}, g^{(m_\nu+1)} \}$,
with $H_\nu \in \{ f^{(m_\nu)} , g^{(m_\nu)} \}$ if $m_\nu$ is exceptional. 
It is then possible to choose 
$j, k \in \{ 1, 2, 3 \}$ with $j < k$ such  that $H_j$ and $H_k$ are both derivatives of $f$, or both derivatives of $g$. 
It follows from (\ref{mjchoose}) that $H_k = H_j^{(q)}$ for some  $q \geq 1$.

Suppose that $m_j$ is exceptional. Then $H_j$ is $f^{(m_j)}$ or $g^{(m_j)}$ and,
because $1 \leq m_j \leq m_k - 2 \leq M-2$,  the choice of $\alpha^*$ and $\beta$ forces $\alpha^* > 0$ 
and $\alpha^* \geq \beta$. Therefore, in this case, (\ref{noto(r)}) yields
$$
\frac14 \, 
T(r_n, H_j'/H_j) \geq \frac{\alpha^* r_n }8   \geq \frac{\beta  r_n }8   \quad \hbox{  for large $n$.}
$$

Thus, whether or not $m_j$ is exceptional, (\ref{choosek}) and (\ref{nra1b1}) or (\ref{nra1b}) give $z_n$ such that (\ref{nrd2}) holds with  $h=H_j$.
Moreover, since $m_k+1 \leq M+1$ and $f, \ldots, f^{(M+3)}$ have finitely many non-real zeros,
combining Lemma \ref{lemnorfam} with (\ref{nra1b1}) or   (\ref{nra1b})
for $m=m_k$ gives (\ref{nrd1}), for some sequence $\rho_n \to 0$. 

Lemma \ref{lemlogder3} now implies that  for large $n$ the set
$$
\left\{ \theta \in  [ \varepsilon ,  \pi - \varepsilon ] : \left| t_n  \frac{H_j'(t_n e^{i \theta } )}{H_j(t_n e^{i \theta } )} \right| < C_1 S_n  \right\} ,
\quad t_n = K^{-1} r_n  , 
$$
has linear measure at least $\pi/2$. On combination with (\ref{choosek}) this yields, as $n \to \infty$, 
\begin{eqnarray*}
 m \left(t_n, H_j/H_j'\right) &\geq& \frac14 \, \log \frac1{S_n} - O(1) = \frac{  \beta r_n }{32} - O(1) = \frac{  \beta K t_n }{32} - O(1) \geq 2 \alpha t_n ,
\end{eqnarray*}
which contradicts (\ref{nrc2b}). 
\hfill$\Box$
\vspace{.1in}

\section{Proof of Theorem \ref{thm1} }

Assume the hypotheses of Theorem \ref{thm1}. Then the results of Section \ref{inter} hold, with $M=9$. 
If $f$ has finite order and finitely many poles then both conclusions of Theorem \ref{thm1} follow from Lemmas \ref{lemfpoly},
\ref{lemfo}  and \ref{lemf'frational3}.
Assume henceforth  that $f$ has infinite order or infinitely many poles and let $Y$ be as in Lemma \ref{lemnomj}.
It will be shown that there exist integers  $m_1, m_2, m_3 \in Y$   satisfying 
(\ref{mjchoose}), contradicting Lemma \ref{lemnomj}.

Suppose first that some $m \in \{ 1, \ldots, 8 \}$ is exceptional.
Then $f^{(m'+1)}$ has infinitely many zeros for $0 \leq m' \neq m$,
by Lemma \ref{lemonenozeros}. 
If $m \leq 3$ then at least one of $a_5$ and $a_6$ is non-zero in (\ref{nra1}),  by Lemma \ref{lemFmFm+1},
as is at least one of $a_8$ and $a_9$: thus $m_1 = m$, while $m_2 \in \{ 5, 6 \}$ and
$m_3 \in \{ 8, 9 \}$. Similarly, if $4 \leq m \leq 6$ then  at least one of $a_1$ and $a_2$ is non-zero, as is at least one of $a_8$ and $a_9$.
Furthermore, if $m \geq 7$ then 
at least one of $a_1$ and $a_2$ is non-zero, as is at least one of $a_4$ and $a_5$.

Suppose finally that $f^{(m+1)}$ has infinitely many zeros, for $m=1, \ldots, 8$. Then Lemma \ref{lemFmFm+1} implies that at least one of $a_2$ and $a_3$ is non-zero
in (\ref{nra1}), as are at least one of $a_5$ and $a_6$ and at least one of $a_8$ and $a_9$. 
\hfill$\Box$
\vspace{.1in}

\section{Proof of Theorem \ref{thm2} }

Assume the hypotheses of Theorem \ref{thm2}. Again the results of Section \ref{inter} hold, this time with $M= 6$, 
and by Lemmas \ref{lemfpoly},
\ref{lemfo}  and \ref{lemf'frational3} it suffices to consider the case where $f^{(m)}$ has infinite order or infinitely many poles, for each $m \geq 0$.
Lemmas \ref{Kmexist} and \ref{lemFmFm+1} imply that if $0 \leq m \leq 7$ and $a_m = 0$ in (\ref{nra1}) then $K_m$ in
(\ref{nrc3})  is  constant and $f^{(m)}$ has finitely many zeros. 

The following argument gives integers  $m_1, m_2, m_3 \in Y$   satisfying 
(\ref{mjchoose}), where $Y$ is as in Lemma \ref{lemnomj}, and so delivers a contradiction.
Suppose first that some $m \in \{ 1, \ldots, 5 \}$ is exceptional, so that $f^{(m'+1)}$ has infinitely many zeros for $0 \leq m' \neq m$, by Lemma \ref{lemonenozeros}.
This implies that if $m \leq 2$ then $a_4 a_6 \neq 0$,  by Lemma \ref{lemFmFm+1}, while if $3 \leq m \leq 4$ then $a_1 a_6 \neq 0$,
and if  $m=5$ then $a_1 a_3 \neq 0$. On the other hand, if  no $m \in \{ 1, \ldots, 5 \}$ is exceptional, then 
$f^{(m+1)}$ has infinitely many zeros, for $m= 1, \ldots, 5$, and $a_2 a_4 a_6 \neq 0$.
\hfill$\Box$
\vspace{.1in}

\section{Proof of Theorem \ref{thmsnr17}}

Let $f$ be a strictly non-real meromorphic function in the plane such that 
all but finitely many zeros and poles of $f$ and $f''$ are real. Write
\begin{equation}
 \label{1a}
 g = \widetilde f, \quad \frac{f'}{f} = \alpha + i \beta, \quad \frac{g'}{g} = \alpha - i \beta , \quad 2 \alpha = \frac{f'}{f} + \frac{g'}{g} ,
\quad 2i \beta = \frac{f'}{f} - \frac{g'}{g} ,
\end{equation}
in which $\alpha$ and $\beta$ are real meromorphic functions. Here  $ \beta $ is not identically zero, since $f/g$ is non-constant,
but $\beta$ has finitely many poles.
Furthermore, all  poles of $\alpha$ are simple, and all but finitely many are real zeros or poles of $f$.
Since $f''/f$ and $g''/g$ have, with finitely many exceptions,  the same zeros and poles there exists a meromorphic function 
$H$ with finitely many zeros and poles such that 
\begin{eqnarray}
\frac{f''}{f}
&=& 
 \alpha' + \alpha^2 - \beta^2 +i( \beta' + 2 \alpha  \beta ) = \frac{H g''}{g}  \nonumber \\
 &=& H 
 \left( \alpha' + \alpha^2 - \beta^2 - i( \beta' + 2 \alpha  \beta ) \right), \quad \widetilde H = \frac1H .
\label{f''frep}
\end{eqnarray}
In view of Lemmas \ref{lem1} and \ref{lem2}, as well as standard properties of the Tsuji characteristic,
\begin{equation}
 \label{f'fbetaest}
T_0(r, f'/f) + T_0(r, g'/g) + T_0(r, \beta) = O( \log r)  \quad \hbox{and} \quad  T(r, \beta ) + T(r, H) = O( r \log r) 
\end{equation}
as $r \to \infty$. 
If $f$ has finite lower order then $\beta$ is a rational function.

Now  $H \equiv 1$ implies that $f''/f$ is real meromorphic and $f'/f$ is a rational function, by  \cite[Theorem 1.3]{LaJdA2013},
and so (\ref{f'festa}) evidently holds: moreover, the same result shows that if, in addition,  $f$ and $f''$ have only real zeros and poles then
$f$ satisfies  (\ref{frepa})(a).

Assume henceforth that $H \not \equiv 1$. Then
rearranging  (\ref{f''frep}) delivers 
\begin{equation}
 \label{4a}
 \alpha' + \alpha^2 - \beta^2 = C( \beta' + 2 \alpha \beta ), \quad C = i \left( \frac{H+1}{H-1} \right) ,
\end{equation}
in which $C$ is a real meromorphic function.

\begin{lem}
 \label{lemalpha0}
If $z_0 \in \C$ is a pole of $\alpha$ but not of $\beta$, and if ${\rm Res} \, (\alpha, z_0)  \neq 1$,  then 
$C(z_0) = \infty $. This holds in particular if $|z_0|$ is large and $z_0$ is a pole or multiple zero of~$f$.
\end{lem}
\textit{Proof.} The residue condition implies that $z_0$ is a 
double pole of $\alpha' + \alpha^2 $, 
and hence a pole of $C$, by (\ref{4a}). The second assertion follows from
(\ref{1a}).
\hfill$\Box$
\vspace{.1in}

Now (\ref{4a}) yields
$$
0= 
\alpha' - C \beta' - C' \beta + \alpha^2 - 2 \alpha C \beta + C^2 \beta^2 + C' \beta - (1+C^2) \beta^2 
$$
and so
\begin{equation}
 \label{5aa}
 0= \gamma' + \gamma^2 + C' \beta - (1+C^2)   \beta^2 , \quad \gamma = \alpha - C \beta . 
\end{equation}

\begin{lem}
 \label{lemHrational}
Assume that $H$ is a rational function in (\ref{f''frep}). Then $f$ 
satisfies (\ref{f'festa}).

If, in addition, $f$ has finite lower order and all  zeros and poles of $f$ and $f''$ are real, then $\beta$, $\gamma$, $\alpha$ and $f'/f$ are all constants,
and $f$ satisfies the first equation of (\ref{frepa}).
\end{lem}
\textit{Proof.} Since $H$ is a rational function, so is $C$. By (\ref{1a}) and Lemma \ref{lemalpha0}, all but finitely many poles
of $\alpha$ are real and simple with residue $1$, and the same is true of $\gamma$ by (\ref{5aa}).
Let $x_0$ be large and positive, and choose $x_1 > x_0$ such that
$\gamma (x_1) \neq \infty$. 
The Riccati equation (\ref{5aa}) may be linearised by writing 
\begin{equation}
 \label{5aaa}
U(x_1) = 1, \quad 
 \frac{U'}{U} = \gamma , \quad 
U'' + (  C' \beta - (1+C^2)  \beta^2 ) U = 0. 
\end{equation}
Then $U$  extends to be analytic in the  half-plane $H_0$ given by ${\rm Re} \, z > x_0$, and $U$ is 
real on $(x_0, \infty)$. 
For $x > x_0$, write $ C'(x) = \rho (x) C(x)$, where $\rho(x)$  is small and real, so that  
$$
1 \geq \frac{ \rho^2 }4  - \left(C \beta - \frac{\rho}2 \right)^2 = \rho  C \beta - C^2 \beta^2  \geq \rho  C \beta - (1+C^2) \beta^2  = C' \beta - (1+C^2) \beta^2 .
$$
Thus the Sturm comparison theorem \cite[p.355]{dureninvite}
applied to $U(x)$ and $V(x) = \sin x$ implies that the number of zeros of $U$ in the interval $[x_0, x]$ is $O(x)$ as 
$x \to + \infty$, and the same is true for the number of poles of $\gamma$, and hence of $\alpha$ and $f'/f$, by (\ref{1a}) and
 (\ref{5aa}). Applying a similar argument on the negative real axis
proves the first estimate of (\ref{f'festa}), and the second follows using (\ref{f'fbetaest}) and Lemma \ref{lem2}. 

Suppose  in addition that
$f$ has finite lower order and all  zeros and poles of $f$ and $f''$ are real. Then $\beta$ is a polynomial in (\ref{1a}), and the rational function
$H$ is free of zeros and poles, 
and so is constant, as is $C$. Moreover, all poles of $\gamma$ are real and simple with residue $1$, so that
$U$ is now a real entire function, with only real zeros,  of finite order by (\ref{5aaa}). 
Furthermore,  $U$ has at most one zero, by the Sturm comparison theorem applied to $U(x)$ and  $V(x) = 1$. Thus 
$\gamma = \alpha - C \beta = U'/U$
has at most one pole, and so is a rational function. 
Hence there exist a polynomial $K = \beta \, \sqrt{1+C^2} \not \equiv 0 $ and a constant $\eta = \pm 1$ such that,  as $z \to \infty$,
(\ref{5aa}) delivers  $\gamma(z) = O( |K(z)|)$ and
\begin{eqnarray*}
 K(z)^2 &=& \gamma(z)^2 + \gamma (z)  \cdot \frac{O(1)}z = \gamma(z)^2 + K(z)  \cdot \frac{O(1)}z,\\
\quad \gamma (z) &=& \eta  K(z) + X(z) = \eta K(z) + \frac{O(1)}z ,\\
0 &=& \eta K'(z) + X'(z) +  2 \eta K(z) X(z) + X(z)^2 = \eta K'(z) +  2 \eta K(z) X(z) +  \frac{O(1)}{z^2} ,
\end{eqnarray*}
as well as 
$$
\frac{U'(z)}{U(z)} + \frac{K'(z)}{2K(z)} = \gamma (z) +  \frac{K'(z)}{2K(z)} = 
\eta K(z) + 
X(z) + \frac{K'(z)}{2K(z)} = \eta K(z) + \frac{O(1)}{z^2 K(z)} .
$$
The argument principle now shows that $U$ and $K$ have no zeros, and hence $K$ and $\beta$ are constant, while 
$\gamma$ is a polynomial and is also constant, as are 
$\alpha$ and $f'/f$. 
\hfill$\Box$
\vspace{.1in}


Assume henceforth that $H$ is transcendental in (\ref{f''frep}). The next lemma follows immediately from (\ref{f'fbetaest}).

\begin{lem}
 \label{lemf''f}
There exist $a\in \R \setminus \{ 0 \}$ and a rational function $T$ with $|T(x)| = 1$ for all $x \in \R$,
such that $H(z) = T(z) e^{iaz}$.  
\end{lem}
\hfill$\Box$
\vspace{.1in}

It may be assumed that $a = 2$ and $T(\infty) = 1$ in Lemma \ref{lemf''f}, so that 
(\ref{4a}) gives
\begin{equation}
 \label{zetadef}
H(z) = e^{2i \zeta(z)}, \quad \zeta (z) = z + \frac{ \log T(z)}{2i} , \quad 
C(z) = i \left( \frac{H(z)+1}{H(z)-1} \right) = \cot \zeta (z),   
\end{equation}
in which the logarithm is the principal branch, while $\zeta (z)$ is analytic near infinity with $\widetilde \zeta = \zeta $ there. 
Thus (\ref{5aa}) becomes 
\begin{equation}
 \label{5a}
 0= \gamma' + \gamma^2 - (1+C^2)  (\beta \zeta' + \beta^2) = 
 \gamma' + \gamma^2 -  (\beta \zeta' + \beta^2) S^2 , \quad S = \frac1{\sin \zeta  } .
\end{equation}


\begin{lem}
 \label{lemalpha4}
Let $x_0$ be large and positive and let $I \subseteq \R \setminus [-x_0, x_0]$ be an open interval containing no poles of $S(z)$.
Then $I$ contains at most one pole of $f'/f$. 
\end{lem}
\textit{Proof.}
Choose $x_1 \in I$ such that $\gamma(x_1) \neq \infty$ and linearise (\ref{5a})  near $x_1$ by writing 
\begin{equation*}
u(x_1) = 1, \quad \frac{u'}{u} = \gamma, \quad 
 u'' + A u=0,  \quad A = - (  \beta \zeta' + \beta^2 )S^2 .
\end{equation*}
Thus $u$ extends to be analytic on a domain containing $I$, and $u$ is real-valued on $I$. Define a zero-free
comparison function $v$ on $I$ by $v(x_1) = 1$ and 
$$
\frac{v'}v =\frac{\zeta' \cot \zeta}2 - \frac{\zeta''}{2 \zeta'}= \frac{\zeta' C}2 - \frac{\zeta''}{2 \zeta'} ,
$$
so that
\begin{eqnarray*}
 \frac{v''}v 
&=& 
\frac{\zeta'' C}2 - \frac{(\zeta')^2 (1+ C^2)}2 - \frac{\zeta'''}{2 \zeta'} +  \frac{(\zeta'')^2}{2 ( \zeta')^2}
+ \frac{(\zeta')^2  C^2}4 - \frac{\zeta'' C}2  + \frac{(\zeta'')^2}{4 ( \zeta')^2}\\
&=& 
- \frac{(\zeta')^2 (1+ C^2)}2 
+ \frac{(\zeta')^2  (1+C^2 -1) }4 - \frac{\zeta'''}{2 \zeta'}  + \frac{3(\zeta'')^2}{4 ( \zeta')^2} \\
&=& 
- \frac{(\zeta')^2 S^2 }{4 } -  \frac{(\zeta')^2   }4 - \frac{\zeta'''}{2 \zeta'}  + \frac{3(\zeta'')^2}{4 ( \zeta')^2} .
\end{eqnarray*}
Since $\zeta'$ is a real rational function with $\zeta'(\infty) = 1$ and  $x_0$ is large, this gives
$$
A = - (\beta \zeta' + \beta^2) S^2 = -\left( \left(\beta+\frac{\zeta'}2\right)^2 - \frac{(\zeta')^2 }4 \right) S^2  \leq  \frac{(\zeta')^2 S^2 }{4 } \leq -\frac{v''}v  
$$
on $I$. 
The Sturm comparison theorem \cite{dureninvite}
now implies that $u$ has at most one zero in $I$, so that $\gamma$ has at most one pole there, as have $\alpha$ and $f'/f$, by (\ref{1a}) and (\ref{5aa}).
\hfill$\Box$
\vspace{.1in}

Since poles of $S$ are poles of $C$ and zeros of $H-1$,  Lemmas \ref{lemalpha0}, \ref{lemf''f} and~\ref{lemalpha4} imply that 
$f$ satisfies  the first estimate of (\ref{f'festa}), from which  the second follows using (\ref{f'fbetaest}) and Lemma \ref{lem2}. 

To complete the proof of Theorem \ref{thmsnr17},  assume henceforth 
that $f$ has finite lower order, all  zeros and poles of $f$ and $f''$ are real and $H$ is transcendental. Then $\beta$ is a polynomial, of degree $d$ say.
Furthermore, $H$ is free of zeros and poles, so that it may be assumed that $H(z) = e^{2i z}$, while $\zeta (z) = z$ and $C(z) = \cot z$. 
Since $\zeta '' \equiv 0$, the next lemma follows from (\ref{4a}), (\ref{5aa}),
Lemma~\ref{lemalpha0} and an 
argument
identical to that in Lemma \ref{lemalpha4}. 

\begin{lem}
 \label{lemalpha1}
(i) Any pole of $f'/f$ in $\C \setminus \pi \Z$ is a simple zero of $f$.\\
(ii) If $z_0 \in \pi \Z$ is a pole of $f'/f$ then ${\rm Res} \, (f'/f, z_0) = 2 \beta (z_0) + 1 $.\\
(iii) If  $n \in \Z$ then $f'/f$ has  in $I_n = (n \pi , (n+1) \pi ) \subseteq \R$ at most one  pole.\\
(iv)
$f$ satisfies 
\begin{equation}
 \label{alphaest2}
N(r, f) + N(r, 1/f) = O( r^{d+1} ) \quad \hbox{as $r \to \infty$.} 
\end{equation}
\end{lem}
\hfill$\Box$
\vspace{.1in}

Now fix $x_1 \in I_0  = (0 ,\pi)$ with $\gamma (x_1) \neq \infty $ and linearise (\ref{5a}) via $u(x_1) = 1$ and $u'/u  = \gamma $, so that $u$ solves
\begin{equation}
 \label{6a*}
u'' + Au = 0, \quad A(z) = - \, \frac{\beta(z) (\beta (z) +1)}{\sin^2 z} .
\end{equation}
Then $u$ extends to be analytic in $\Omega = \C \setminus \{ n \pi - it: \, n \in \Z, \,  t \in [0, + \infty)  \}$, with 
$u$ real on $I_0$. 


\begin{lem}
 \label{lemgamma0}
Let $0 < \varepsilon < \pi /4$ and denote by $E_0(z)$ any term which satisfies $\log^+ |E_0(z)| = o(|z|)$ 
as $z \to \infty$ with $\varepsilon < \arg z < \pi - \varepsilon $. 
Then there exists   a polynomial $P \not \equiv 0$ of degree at most $1$ such that
\begin{equation}
 \label{gammaest1}
\frac{u''(z)}{u(z)} = E_0(z) e^{2  i z} , \quad 
u(z) = P(z) + E_0(z) e^{2iz} , \quad 
\gamma (z) = \frac{P'(z)}{P(z)} + E_0(z) e^{2iz} . 
\end{equation}
\end{lem}
\textit{Proof.} The first estimate follows from (\ref{6a*}) and the remaining two are proved by the method of Gronwall's lemma, exactly as in 
\cite[Lemma 4.3]{LaJdA2013}.
\hfill$\Box$
\vspace{.1in}

\begin{lem}
 \label{lemorderf}
The order of $f$ is at most $d+1$. 
\end{lem}
\textit{Proof.} (\ref{alphaest2}) makes it possible to write $f = \Pi e^Q$ where $\Pi$ is a meromorphic function with real zeros and poles and order
at most $d+1$, while $Q$ must be a polynomial. It follows from (\ref{1a}), (\ref{5aa}), (\ref{gammaest1}) and standard estimates for logarithmic derivatives
that, as $z \to \infty$ with $\varepsilon < \arg z < \pi - \varepsilon $,
$$
Q'(z) = \frac{f'(z)}{f(z)} - \frac{\Pi'(z)}{\Pi(z)} = \gamma (z) + (\cot z + i ) \beta(z) - \frac{\Pi'(z)}{\Pi(z)} = O( |z|^{d + 1/2} ) ,
$$
so that $Q$ has degree at most $d+1$. 
\hfill$\Box$
\vspace{.1in}

\begin{lem}
 \label{lemgamma1}
If the degree $d $ of $\beta$ is positive then, as $x \to + \infty$ with $x \in \R$, 
\begin{equation*}
\left| \left( \frac{f'}{f} \right)' (x+i) \right| + 
| \alpha'( x + i)| +  | \gamma'( x + i)| = o( | (x+i) \beta (x+i)| ) = o( | \beta (x+i)|^2 ) .
\end{equation*}
\end{lem}
\textit{Proof.} 
It suffices by (\ref{1a}) and (\ref{5aa}) to prove that $(f'/f)'(x+i) = o( | (x+i) \beta (x+i)| )$. 
Let $x \in (0, + \infty )$ be large, set $w = x+i$ and take  $R \in [ 2 |w|, 2|w| + 1]$ such that $f(z) \neq 0, \infty$ on $|z|=R$.
Denote by $a_j$ the zeros and poles of $f$ in $|z| < R$, repeated according to multiplicity. 
Applying the twice differentiated Poisson-Jensen formula \cite[(1.17)]{Hay2}  to $f$ in the disc $|z| < R $  gives 
\begin{eqnarray*}
 \left| \left( \frac{f'}{f} \right)' (w) \right| &\leq& 
\frac2{\pi} \int_0^{2 \pi} \frac{ R | \log |f(Re^{it})| | }{| Re^{it} - w|^3 } \, dt + 
\sum \left( \frac1{|a_j-w|^2} + \frac{|a_j|^2}{|R^2 - \overline{a_j} w|^2 } \right) ,
\end{eqnarray*}
in which $| Re^{it} - w| \geq R/2$, while $|R^2 - \overline{a_j} w| \geq (1/2) R^2 $ and $|a_j-w| \geq 1$.
Lemma~\ref{lemalpha1} implies that the number of distinct zeros and poles of $f$ in the interval
$[x-  R/ \log R , x+  R/ \log R ]$ is $O( R/\log R)$, and that each of these
has multiplicity at most $4 M(R, \beta) $.
It now follows from Lemma~\ref{lemorderf} that 
\begin{eqnarray*}
 \left| \left( \frac{f'}{f} \right)' (w) \right| &\leq& \frac{32}{R^2} (m(R, f) + m(R, 1/f) ) + 
 O \left( \frac{R  \, M(R, \beta) }{ \log R} \right)  + \\
& & +  (n(R,f) + n(R, 1/f)) \left( \frac{(\log R)^2}{R^2} + \frac4{R^2} \right) \\
&\leq&  O( R^d ) + O \left( \frac{R  \, M(R, \beta) }{ \log R} \right)  = o( R \, M(R, \beta ) ) = o(   | w \beta (w)| ).
\end{eqnarray*}
\hfill$\Box$
\vspace{.1in}

\begin{lem}
\label{lembetaconst}
The polynomial $\beta$ has degree $d=0$ and, without loss of generality, there exists a real meromorphic function $W$ on $\C$ of order at most $1$
such that  
\begin{equation}
 \label{Wdefa}
 f(z) = W(z)e^{i\beta z} , \quad \frac{W'}{W} = \alpha = \gamma + \beta C = \frac{u'}{u} + \beta C .
\end{equation}
\end{lem}
\textit{Proof.} Assume that $\beta$ has positive degree $d $ and let $\varepsilon $ be small and positive. 
The equations 
(\ref{1a}) and (\ref{5aa}) and 
the fact that $f$ has finite order give $M_2 > 0$ and arbitrarily large positive $R$ with 
$\gamma(z) = O\left( R^{M_2} \right)$ on $|z| = R$. 
Now 
Lemma \ref{lemgamma0} shows that 
$$
\frac{ ( \gamma (z) - P'(z) /P(z) ) \sin z }{\beta (z) } \to 0
$$
as $z \to \infty$ with $\arg z = 2 \varepsilon $, whereas (\ref{6a*}) and Lemma \ref{lemgamma1}  imply that 
$$
\gamma (x+i)  \sim \pm \, \frac{\beta(x+i)}{\sin (x+i) }  , \quad 
\frac{ ( \gamma (x+i) - P'(x+i)/P(x+i)  ) \sin (x+i) }{\beta (x+i) } \to \pm 1 ,
$$
as $x \to + \infty$ with $x \in \R$. 
Since 
$\gamma$ has only real poles,
this contradicts the Phragm\'en-Lindel\"of principle. 
The remaining assertions follow from (\ref{1a}), (\ref{5aa}) and Lemma \ref{lemorderf}.
\hfill$\Box$
\vspace{.1in}

\begin{lem}
 \label{lembeta1}
 If $u(z)$ and $u(z+\pi)$ are linearly dependent on $\Omega$  then $f$ satisfies (\ref{frepa}). 
\end{lem}
\textit{Proof.} The hypotheses imply that $\gamma = u'/u$ has period $\pi$ 
and so have
the sequences of poles and zeros of $f$, by (\ref{1a}) and (\ref{5aa}).
Thus, by  Lemma~\ref{lemalpha1},
either $f$ has in each interval $I_n = (n \pi, (n+1) \pi )$, $n \in \Z$,  
exactly one simple zero and no poles, or $f$ has no zeros and poles in the $I_n$.  Moreover, the residue of $f'/f$ at each zero of $\sin z$ is 
a fixed integer $m$, possibly $0$. It follows that $f$ has a representation 
\begin{equation}
 \label{flastrep}
f(z) =  (e^{2iz} - 1)^L (e^{2iz} - E) e^{pz+q} , \quad L \in \Z,  \quad E, p, q \in \C, \quad |E|=1,
\end{equation}
in which $E = 1$ is not excluded. 
This implies in view of (\ref{1a}) and (\ref{5aa}) that, as $z \to \infty$ in $\varepsilon < \arg z < \pi - \varepsilon $, 
$$
\frac{f'(z)}{f(z)} = p  + o(1), \quad \alpha (z) = p - i \beta + o(1),
\quad \gamma (z) = \alpha (z) - \beta \cot z = p + o(1) ,
$$
so that $p= 0$ by Lemma \ref{lemgamma0}. Now $L \neq -1$ in (\ref{flastrep}), since $f$ is strictly non-real, and $f$ is determined by applying Lemma \ref{lemf1}
to $F(z) = e^{-q} f(z/2i)$. 
\hfill$\Box$
\vspace{.1in}

Assume henceforth that  $u(z)$ and $u(z+\pi)$ are linearly independent solutions on~$\Omega$ of (\ref{6a*}). 
The proof of Theorem \ref{thmsnr17} will be completed by first considering certain values of $\beta$ with $|\beta|$  small, 
following which the remaining possibilities for $\beta$ will be disposed of together.

\begin{lem}
 \label{lembeta3}
 If $ \beta \in \{ -2, -1, 1 \}$ then $f$ satisfies (\ref{frepa}).
\end{lem}
\textit{Proof.} Suppose first that $\beta = -1$: then (\ref{6a*}) shows that  $u''=0$.
By (\ref{Wdefa}) and the fact that $u(z)$ and $u(z+\pi)$ are linearly independent, there exists a polynomial $T_1$, of 
degree $1$, such that
$$
\gamma = \frac{u'}u = \frac{T_1'}{T_1} , \quad 
f(z) = W(z) e^{-iz} = \frac{T_1(z)}{e^{2iz} - 1} .
$$
Now (\ref{f''frep}) and (\ref{5aa}) lead to
\begin{eqnarray*}
 \frac{f''}{f} &=& \gamma' +(1+C^2)  + \gamma^2 -2C \gamma + C^2 - 1 -  2i (\gamma - C) \\
&=&  -2C \gamma + 2 C^2  -  2i (\gamma - C) = 2(C+i)(C - \gamma) . 
\end{eqnarray*}
Since $f''$ has only real zeros, all zeros of $C - \gamma$ must be real. Thus 
the zero of $T_1$ belongs to $ \pi \Z$; if this is not the case then Lemma \ref{lemtan} gives a   non-real zero $z^*$ of
$\tan z - T_1(z)/T_1'(z)$, with $\tan (z^*) \neq 0, \infty $ and so $T_1(z^*) \neq 0$, 
which implies that $z^*$  is a non-real zero of $C-\gamma$, a contradiction.
It follows that $f$ is given by (\ref{frepa})(c). 

Now suppose that $ \beta \in \{ -2,  1 \}$. Then $\beta (\beta + 1) = 2$  and 
(\ref{6a*}) solves explicitly to give
$A_1, B_1 \in \C$ with 
$$
u(z) = A_1 \cot z + B_1 ( 1 - z \cot z ) ,
$$
in which $B_1 \neq 0$ since $u(z)$ and $u(z+\pi)$ are linearly independent. Hence there exists 
a polynomial $T_1$ of degree $1$ such that,
in view of (\ref{Wdefa}), 
\begin{equation}
 \label{fbetarep}
f(z) = ( T_1'(z) - T_1(z) \cot z) (\sin z)^\beta  e^{i \beta z} .
\end{equation}
If  $\beta = 1$ this gives (\ref{frepa})(b), and again the zero of $T_1$ must belong to $\pi \Z$ by Lemma \ref{lemtan}.
 
Assume now that $\beta = -2$. Then 
Lemma \ref{lemalpha1} implies that $f$ has no multiple zeros. Suppose that $x_0 \in \R$ is a  simple zero of $f$, and so a simple pole with residue $1$ of
the real meromorphic function $\alpha$. Then there exists $D_0 \in \R$ such that, as $z \to x_0$, 
$$
\frac{f'(z)}{f(z)} = \alpha (z) + i \beta = \frac1{z-x_0} + D_0 - 2i + O (|z-x_0|) , \quad 
\frac{f''(z)}{f(z)} = \frac{2(D_0-2i)}{z-x_0} + O(1). 
$$
This shows that $x_0$ is a pole of $f''/f$, and so not a zero of $f''$.
Thus every  zero of $f''$ must be a real zero of $f''/f$ and so of $\alpha$, by
(\ref{f''frep}). But (\ref{fbetarep}) leads to
$$
\quad \alpha = \frac{f'}{f} - i \beta = \frac{-T_1'C + T_1(1+C^2)}{T_1'-T_1C} -2C = \frac{- 3T_1'C + T_1 + 3 T_1 C^2 }{T_1'-T_1C} .
$$
Hence if $|z|$ is large and $\alpha(z)=0$ then $C \neq \infty$ and
$3C^2 + 1 = o(1) C$, so that $C$ is non-real and so is $z$.
Therefore $f''$ has finitely many zeros and, by the main result of \cite{Lanew}, $f$ has finitely many poles,  contradicting (\ref{fbetarep}). 
\hfill$\Box$
\vspace{.1in}

\begin{lem}
 \label{lembeta2}
Let $n \in \Z$. Then  near $n \pi$ there exist linearly independent local solutions $u_1, u_2$  of  (\ref{6a*}) of form 
\begin{equation}
\label{u1u2def}
 u_1(z) = (z-n \pi)^{-\beta} h_1(z), \quad u_2(z) = (z-n \pi)^{\beta+1} h_2(z) , \quad 
 h_j(z) = 1 + \sum_{k=1}^\infty a_{j,k} (z - n \pi)^k , 
\end{equation}
in which the $h_j$ are analytic on $|z-n \pi| < \pi$ and
the coefficients $a_{j,k}$ are independent of $n$. 
Moreover,  $2 \beta + 1 $ is an integer, and $\beta \neq \pm 1/2$ and $\beta \neq -3/2$. 
Finally, if $u_3, u_4$ are non-trivial solutions on $\Omega$ of  (\ref{6a*}), then $u_3^2, u_4^2 $ and $u_3/u_4$ 
 all extend to be meromorphic in the plane. 
\end{lem}
\textit{Proof.} Choose some $n \in \Z$ and observe first that, near the  regular singular point 
$n \pi$, there exists $ \delta  \in \{ -\beta, \beta+1 \}$ such that (\ref{6a*}) has a solution of form 
\begin{equation}
 \label{U1def}
U_1(z) = (z-n \pi)^{\delta} H_1(z),  \quad H_1(z) =  1+ \sum_{k=1}^\infty b_k (z - n \pi)^k , 
\end{equation}
with $H_1$ analytic on $|z-n \pi| < \pi$.  Since $U_1(z+\pi)$ solves  (\ref{6a*}), for $z$ near $(n-1)\pi$, such a solution exists for any $n$,
with the same choice of $b_k$. 
To obtain a further solution near $n \pi$ write
\begin{eqnarray}
 \label{U2def}
U_2(z) &=& U_1(z) \int U_1(z)^{-2} \, dz 
= U_1(z) \int (z-n \pi )^{-2 \delta } (1 -2 b_1 (z-n \pi)  + \ldots ) \, dz \nonumber \\
&=& U_1(z) \left( c_1 \log (z-n \pi) + (z-n \pi )^{1 -2 \delta } (d_0 + d_1(z- n \pi) + \ldots ) \right) ,
\end{eqnarray}
in which the series $\sum_{k=0}^\infty d_k (z- n \pi)^k $ is obtained by formal integration but has positive radius of convergence. 
Suppose first that $c_1 \neq 0$.
Then $-2 \delta \in \Z$, and so  $U_1^{-2} = (U_2/U_1)'$ has a  meromorphic extension to a neighbourhood of $n \pi$,
as have $\gamma = u'/u$ and $U_1'/U_1$.
Write the solution $u$ of  (\ref{6a*}) locally in the form $u = \alpha_1 U_1 + \alpha_2 U_2$ near $n \pi$, with the $\alpha_j \in \C$.
Then the logarithmic derivative of $\alpha_1 + \alpha_2 U_2/U_1$  extends  meromorphically to
a neighbourhood of $n \pi$,  and if $\alpha_2 \neq 0$  so does $U_2/U_1$,  a contradiction. 
Hence $u$ must locally be a constant multiple of $U_1$ only, 
so that $u(z)$ and $u(z+\pi)$ are linearly dependent, contrary to assumption. 

Thus a logarithm cannot arise in (\ref{U2def}), 
which forces $\beta \neq -1/2$ and $- \beta \neq \beta + 1$, and there exist 
local solutions $u_1, u_2$ as in (\ref{u1u2def}), obtained via (\ref{U1def}) and (\ref{U2def}), 
with 
the coefficients $a_{j,k}$ independent of $n$.
It follows that near $n \pi$ the meromorphic function $W$ in  (\ref{Wdefa}) is a linear
combination of 
$$
v_1(z) = k_1(z), \quad v_2(z) = (z-n \pi)^{2\beta+1} k_2(z) , 
$$
where the $k_j$ are analytic on $|z-n \pi| < \pi$, with $k_j(n \pi ) \neq 0$. But then, if $2 \beta + 1 \not \in \Z$,
it must be the case that $W$ is a constant multiple of $v_1$ only, 
so that $u(z)$ and $u(z+\pi)$ are linearly dependent,  again contrary to assumption.

Next, suppose that $\beta = 1/2$ or $\beta = -3/2$. Then one of $-\beta$ and $\beta + 1$ is $-1/2$ and  by (\ref{u1u2def})
there exists,   near $0$, a solution  of  (\ref{6a*}) of form
$U_3(z) = z^{-1/2} (1 + e_1 z + e_2 z^2 + \ldots )$,
so that 
\begin{eqnarray*} 
 \beta ( \beta + 1)  U_3(z) &=& \frac34 \, z^{-1/2} (1 + e_1 z + e_2 z^2 + \ldots )  = U_3''(z) \sin^2 z \\
&=& \left( \frac34 \, z^{-5/2} - \frac14 \, e_1 z^{-3/2} + \frac34 \, e_2 z^{-1/2} + \ldots \right) \left( z^2 - \frac{z^4}3 + \ldots \right)  \\
&=& \frac34 \,  z^{-1/2} - \frac14 \, e_1 z^{1/2} + z^{3/2} \left( \frac34 \, e_2 - \frac14  \right)  + \ldots .
\end{eqnarray*}
Comparing the coefficients of $z^{3/2}$ yields a contradiction. 

To complete the proof observe that, because $2 \beta + 1  \in \Z$, the  $u_j$ in (\ref{u1u2def}) are such that 
$u_1^2, u_2^2 $, $u_1 u_2$  and $u_1/u_2$  extend to be meromorphic on a neighbourhood of $n \pi \in \pi \Z$.  

\hfill$\Box$
\vspace{.1in}



In view of Lemmas \ref{lembeta3} and \ref{lembeta2}, as well as the fact that $f$ is strictly non-real,
it remains only to consider the case where 
$2 \beta + 1 \in \Z$ but 
\begin{equation}
 \label{betanot}
\beta \not \in \{ -2, -3/2, -1, -1/2, 0, 1/2,  1 \}, \quad \beta + 1 \not \in \{ -1, -1/2, 0, 1/2, 1, 3/2,  2 \} .
\end{equation}

\hfill$\Box$
\vspace{.1in}


\begin{lem}
 \label{lembeta4}
If  $n \in \Z$ then $u^2$ has at $n \pi$ a zero or pole of multiplicity at least $3$.
 
Furthermore, there exist infinitely many $n \in \Z$ such that $n \pi $ is a pole of $u^2$, 
and infinitely many $n \in \Z$ such that $n \pi $ is a zero of $u^2$.
\end{lem}
\textit{Proof.} The   equation  (\ref{6a*}) has local solutions $u_j$ as in
(\ref{u1u2def}), in which $-\beta$ and $\beta+1$ have opposite signs, and 
$u^2$ has a zero or pole at $n \pi$ of  multiplicity $2 | \beta |\geq 3 $ or $2 | \beta + 1| \geq 3$, by (\ref{betanot}).

To prove the last assertion, assume that $u^2$ has a pole at all but finitely many $n \pi$, $n \in \Z$, or that
$u^2$ has a zero at all but finitely many of these points. In the first case set $V = u^2$, and in the second set $V = u^{-2}$. 
Then $V$ satisfies, as $r \to \infty$, 
\begin{equation}
 \label{polest1}
\frac{6r}\pi - O(1) \leq n(r, V), \quad \frac{6r}\pi - O( \log r) \leq N(r, V)  .
\end{equation}
On the other hand, Lemma \ref{lemalpha1} shows that if $n \in \Z$ then in the interval $I_n = (n \pi, (n+1)\pi)$ the function
$f'/f$ has at most one pole, and any such pole
has residue $1$. The same is true 
of $\alpha$ and $\gamma = u'/u$, by (\ref{1a}) and (\ref{5aa}),  and so $u^2$ has no poles and 
at most two zeros in $I_n$. This implies that,  as $r \to \infty$, by (\ref{polest1}) and Jensen's formula,
\begin{equation}
\label{zeroest1}
 N(r, 1/V ) \leq \frac{4r}\pi + O( \log r) \leq \left( \frac23 + o(1) \right) N(r, V) \leq \frac34 \, T(r, 1/V) .
\end{equation}

Since $f$ has finite order, 
applying \cite[Lemma 4]{Hay3} 
gives $C_0 > 1$ and a set $E_1 \subseteq [1, \infty)$, of
positive lower logarithmic density, such that  $T(2r, 1/V) \leq C_0 T(r, 1/V)$ for all $r \in E_1$. 
Choose a positive constant $\varepsilon$, so small that
$$
88 C_0     \varepsilon    \left( 1 + \log^+ \frac1{4 \varepsilon } \right)  <  \frac1{16} .
$$
Then  Lemma \ref{lemgamma0}, the fact that $u^2$ is real meromorphic and
an inequality of Edrei and Fuchs 
\cite[ p.322]{EF2} together deliver, for large $r \in E_1$,
\begin{eqnarray*}
m(r, 1/V) &\leq& O( \log r) + 11 \left( \frac{2r}{2r-r} \right) \, 4 \varepsilon   \, \left( 1 + \log^+ \frac1{4 \varepsilon } \right)  T(2r, 1/V)  \\
&\leq& O( \log r) + 
88 C_0   \varepsilon    \left( 1 + \log^+ \frac1{4 \varepsilon } \right) T(r, 1/V) \leq \frac18 \, T(r, 1/V) , 
\end{eqnarray*}
which contradicts (\ref{zeroest1}).
\hfill$\Box$
\vspace{.1in}

\begin{lem}
\label{lemperiodic}
The function 
$$
G(z) = \frac{u(z+\pi)-u(z)}\pi 
$$
is a non-trivial solution of   (\ref{6a*}) with period $\pi$ on  $\Omega $.
\end{lem}
\textit{Proof.} Lemma \ref{lemgamma0} shows that $u(z)$ is asymptotic to a polynomial $P \not \equiv 0$ of degree at most~$1$ as $z \to \infty$ 
in $\varepsilon < \arg z < \pi - \varepsilon $. The Wronskian $W_u$ of $u(z)$ and $u(z+\pi)$ is constant,  by Abel's identity and (\ref{6a*}). If $P$ is constant
then $W_u$ 
tends to $0$ in a sector and so must vanish identically, forcing $u(z)$ and $u(z+\pi)$ to be linearly dependent, contrary to assumption.

Thus $P$ must be non-constant, and $G(z)$ and $G(z+\pi)$  both solve   (\ref{6a*}) and are asymptotic to the same non-zero constant as $z \to \infty$ 
in $\varepsilon < \arg z < \pi - \varepsilon $. The argument of the previous paragraph now shows that
$G(z)$ and $G(z+\pi)$ are linearly dependent and must be equal. 
\hfill$\Box$
\vspace{.1in}

It is now possible to write  
$$
u(z) = z G(z) + K(z) , \quad \frac{K(z)}{G(z)} = \frac{u(z)}{G(z)} -z  ,
$$
where $K$ also has period $\pi$ on $\Omega$. Moreover, $G^2$ and $K/G$ are meromorphic in the plane,  by Lemma \ref{lembeta2}, and have period $\pi$. 
Lemma \ref{lembeta4} implies that $G^2$ has at least one pole in $ \pi \Z$, 
and so a pole at every point of $\pi \Z$.
If $n \in \Z$ and $n \pi$ is not a pole of $u^2$ then,
as $z \to n \pi$ with $z \in \Omega$,
$$
u(z) = z G(z) + K(z) = O(1), \quad 
\frac{K(z)}{G(z)}  =  \frac{u(z)}{G(z)} - z \to - n \pi ,
$$
which cannot hold for more than one such $n \pi$, since $K/G$ is periodic. 
Thus $u^2$ has a pole  at all but at most one $n \pi \in \pi \Z $,  contradicting Lemma \ref{lembeta4}. 
\hfill$\Box$

\noindent
School of Mathematical Sciences, University of Nottingham, NG7 2RD, UK.\\
james.langley@nottingham.ac.uk

\end{document}